\documentclass{amsart}
\usepackage{amssymb}   
\usepackage{amsmath}

\newtheorem{Theorem}{Theorem}[section]
\newtheorem{Lemma}[Theorem]{Lemma}
\newtheorem{Proposition}[Theorem]{Proposition}
\newtheorem{Corollary}[Theorem]{Corollary}

\newtheorem{defn}[Theorem]{Definition}

\newtheorem{Remark}[Theorem]{Remark}
\newtheorem{remark}[Theorem]{Remark}
\newtheorem{Claim}[Theorem]{Claim}

\newtheorem{example}[Theorem]{Example}
\numberwithin{equation}{section}

\newcommand{\mult}{mult}
\newcommand{\cal}{\mathcal}
\renewcommand{\O}{{\mathcal O}}

\renewcommand{\L}{{\mathcal L}}
\newcommand{\Proj}{{\mathbb P}}

\newcommand{\FF}{{\mathbb F}}

\newcommand{\Z}{{\mathbb Z}}

\newcommand{\C}{{\mathbb C}}
\newcommand{\p}{{\mathbb P}}

\newcommand{\codim}{\operatorname{codim}}
\newcommand{\red}{\operatorname{red}}

\renewcommand{\mult}{\operatorname{mult}}

\newcommand{\Deltaprimo}{\operatorname{\Delta'}}
\newcommand{\map}{\dasharrow}

\begin{document}

\title[Minimal secant degree]{Varieties with minimal secant degree and
linear systems of maximal dimension on surfaces}

\author{Ciro Ciliberto}
\address{Dipartimento di Matematica, Universit\'a di Roma Tor Vergata,
Via della Ricerca Scientifica,
 00133 Roma, Italia}
\email{cilibert@axp.mat.uniroma2.it}

\author{Francesco Russo}
\address{Departamento de Matematica, Universidade Federal
de Pernambuco, Cidade Universitaria, 50670-901 Recife--PE, Brasil}
\email{frusso@dmat.ufpe.br}

\subjclass{Primary 14N05; Secondary 14C20, 14M20}

\keywords{Higher secant varieties, tangential projections, linear
systems}

\maketitle

\section*{Introduction}

In this paper, in which we work over the field of complex numbers,
we touch, as the title suggests, two different themes, i.e. secant
varieties and linear systems, and we try to indicate some new,
rich, and to us unexpected, set of relations between them.\par

Let $X\subseteq\p^r$ be a reduced, irreducible, projective
variety. Basic geometric objects related to $X$ are its secant
varieties $S^k(X)$, i.e. the varieties described by all projective
subspaces $\p^k$ of $\p^r$ which are $(k+1)$--secant to $X$ (see
\S \ref {subsecsec} for a formal definition: in \S \ref {notation}
we collected all the notation and a bunch of useful preliminaries
which we use in the paper). The presence of secant varieties in
the study of projective varieties is ubiquitous, since a great
deal of projective geometric properties of a variety is encoded in
the behaviour of its secant varieties. However the importance of
secant varieties is not restricted to algebraic geometry only.
Indeed, different important problems which arise in various fields
of mathematics can be usefully translated in terms of secant
varieties. Among these it is perhaps the case to mention
polynomial interpolation problems, rank tensor computations and
canonical forms, expressions of polynomials as sums of powers and
Waring type problems, algebraic statistics etc. (see, for
instance,  \cite {C}, \cite {CGG}, \cite {GSS}, \cite {IK}).
\par

Going back to projective algebraic geometry, let us mention the
first basic example of a property of a variety which is reflected
in properties of a secant variety: it is well known, indeed, that
a smooth variety $X\subseteq \p^r$ can be projected isomorphically
to $\p^{r-m}$, with $m>0$, if and only if its first secant variety
$S(X):=S^1(X)$ has codimension at least $m$ in $\p^r$. Furthermore
one can ask {\it how singular} a general projection of $X$ to
$\p^{r-m-1}$ from a general $\p^m$ is, if $m$ is exactly the
codimension of $S(X)$ in $\p^r$. One moment of reflection shows
that a basic step in answering this question is to know in how
many points $S(X)$ intersects a general $\p^m$ in $\p^r$, i.e. one
has to know what is {\it degree} of $S(X)$. A related, more
difficult problem, is to understand what is the structure of the
cone of secant lines to $X$ passing through a general point in
$S(X)$, a classical question considered by various authors even in
very recent times (see, for instance, \cite {LR}). Of course
similar problems arise in relation with higher secant varieties
$S^k(X)$ as well and lead to the important questions of
understanding what is the {\it dimension} and the {\it degree} of
$S^k(X)$ for any $k\geq 1$. \par

As well known, if $X$ has dimension $n$, there is a basic upper
bound for the dimension of $S^k(X)$ which is provided by a naive
count of parameters (see (\ref {defect}) below). As often happens
in many similar situations in algebraic geometry, one expects that
{\it most} varieties achieve this upper bound, and that it should
be possible to classify all the others, the so--called $k$--{\it
defective varieties}, namely the ones for which the dimension of
$S^k(X)$ is smaller than the expected. Unfortunately this
viewpoint, which is in principle correct, is in practice quite
hard to be successfully pursued. Indeed, while there are no
defective curves and the classification of defective surfaces,
though not at all trivial, is however classical (see \cite {Se},
\cite {Terr2} and \cite {CC} for a modern reference), the
classification of defective threefolds is quite intricate and has
only recently been completed (see \cite {CC3}) after the classical
work of Scorza \cite {Scorza} on $1$--defective threefolds (see
also \cite {CC2}). As for higher dimensional defective varieties,
no complete classification result is available, though a number of
beautiful theorems concerning some special classes of defective
varieties is available (see \cite {Zak}).
\par

One of the objectives of the present paper is to address the other
question we indicated above, i.e. the one concerning the
determination of the degree of secant varieties. This question,
though important, has never been systematically investigated in
general, neither in the past, nor in more recent times, exceptions
being, for instance, the paper \cite {CJ} for the case of curves
(see also \cite{Zak2}), and the computation of the degree of
secant varieties to varieties of some particular classes, like one
does in  \cite {Room} (see also \S \ref{ex} below).\par

Of course, given any variety $X\subseteq \p^r$, one has a famous,
classical lower bound for the degree of $X$ (see (\ref{mindeg})
below), which says that the degree in question is bounded below by
the codimension of $X$ plus one. This bound is sharp, and the
varieties achieving it, the so--called {\it varieties of minimal
degree}, are completely classified, in particular they turn out to
be rational (see  \cite {EH}). The aforementioned bound of course
applies to the secant varieties of $X$ too, but, according to the
classification of varieties of minimal degree, one immediately
sees that it is never sharp in this case. Thus the question arises
to give a {\it sharp} lower bound for the degree of $S^k(X)$. This
is the problem that we solve in \S \ref {minsecdeg}, where our
main result, i.e. Theorem \ref {maxbound}, is the bound
(\ref{goodbound}) for the degree of $S^k(X)$. Moreover we prove a
similar bound (\ref {goodbound2}) for the multiplicity of $S^k(X)$
at a general point of $X$. One of the main steps in the proof of
Theorem \ref {maxbound} is the result in \S \ref {tangentcones},
namely Theorem \ref {conotang}, in which we give relevant
informations about the tangent cone to $S^k(X)$ at the general
point of $S^l(X)$, where $l<k$. This can be seen as a wide
generalization of the famous Terracini's Lemma (see Theorem
\ref{terracini} below), which describes the general tangent space
to $S^k(X)$.\par

The lower bound  (\ref{goodbound}) for the degree of $S^k(X)$ is a
generalization of the classical lower bound (\ref{mindeg}) for the
degree of any variety, and, as well as the latter, it is sharp.
Actually, in Theorem \ref {maxbound} we also show that varieties
$X$ such that $S^k(X)$ has the minimum possible degree, called
varieties with {\it minimal $k$--secant degree} or
$\mathcal{M}^{k}$--{\it varieties} (see Definition \ref {msdeg}),
enjoy important properties like: general $m$--internal projections
$X^m$ of $X$, i.e. projections of $X$ from $m$ general points on
it, are also of {\it minimal $k$--secant degree}, {general
$m$--tangential projections} $X_m$ of $X$, i.e. projections of $X$
from $m\leq k$ general tangent spaces, are of {\it minimal
$(k-m)$--secant degree}, in particular, for $k=m$, projections
$X_k$ of $X$ from $k$ general tangent spaces are of {\it minimal
degree}, hence they are rational. Since we know very well
varieties of minimal degree, and a general $k$--{\it tangential
projection} $X_k$ of $X$ is one of them, a natural question, at
this point, arises: what is the structure of the projection $X\map
X_k$? The interesting answer is that, if $X$ is not $k$--defective
then the map in question is generically finite and its degree is
bounded above by $\mu_k(X)$ which, by definition, is the number of
$(k+1)$--secant $\p^k$ to $X$ passing through the general point of
$S^k(X)$. In particular, if $X$ is not $k$--defective, if $S^k(X)$
has minimal degree and $\mu_k(X)=1$, then  $X$, as well as $X_k$,
is rational. The main ingredient for the proof of the bound on the
degree of the $k$--tangential projection $X\map X_k$ is proved in
\S \ref {degproj} (see Theorem \ref {boh}), where we exploit and
generalize the technique, introduced in \cite {CMR}, of
degeneration of projections, based on a beautiful idea of
Franchetta (see \cite {Fr1},\cite{Fr2}). \par

 Notice that the condition $\mu_k(X)=1$ is rather mild, i.e. one
expects that {\it most} non $k$--defective varieties $X\subset
\p^r$ enjoy this property if $S^k(X)\subsetneq\p^r$ (see \S
\ref{subsecwd}, in particular Proposition \ref{wd} for a
sufficient condition for this to happen). The varieties $X$, not
$k$--defective, such that $S^k(X)$ has minimal degree and
$\mu_k(X)=1$ are called $\mathcal{MA}^{k+1}_{k-1}$--variety or
$\mathcal{OA}^{k+1}_{k-1}$--variety according to whether $S^k(X)$
is strictly contained in $\p^r$ or not (see Definition
\ref{msdeg}), e.g. $X$ is an $\mathcal{OA}^{k+1}_{k-1}$--variety
if and only if $S^k(X)=\p^r$, $r=(k+1)n+k$ and there is only one
$(k+1)$--secant $\p^k$ to $X$ passing through the general point of
$\p^r$, i.e. the general projection $X'$ of $X$ to $\p^{r-1}$
acquires a new $(k+1)$--secant $\p^{k-1}$ that $X$ did not use to
have. This was classically called an {\it apparent $(k+1)$--secant
$\p^{k-1}$} of $X$. It should be mentioned, at this point, the
pioneering work of Bronowski on this subject: in his inspiring,
but unfortunately very obscure, paper \cite {Br1} he essentially
states that the map $X\map X_k$ is birational if and only if $X$
is either an $\mathcal{MA}^{k+1}_{k-1}$--variety or
$\mathcal{OA}^{k+1}_{k-1}$--variety. As we said, one implication
has been proved by us, the other is open in general, and we call
it the $k$--th Bronowski's conjecture (see Remark
\ref{Bronowski}). The results of the present paper imply that
Bronowski's conjecture holds for smooth surfaces (see Corollary
\ref {Bron}), whereas the main theorem of \cite {CMR} implies that
Bronowski's conjecture holds for smooth threefolds in $\p^7$ if
$k=1$. It would be extremely nice to shed some light on the
validity of this conjecture in general, since, according to
Bronowski, this would make the study and the classification of
$\mathcal{MA}^{k+1}_{k-1}$ and
$\mathcal{OA}^{k+1}_{k-1}$--varieties easier. \par

The existence of $\mathcal{M}^{k}$, $\mathcal{MO}^{k+1}_{k-1}$,
and $\mathcal{MA}^{k+1}_{k-1}$--varieties, and therefore the
sharpness of the bound proved in Theorem \ref{maxbound}, is showed
in \S \ref{ex}, where several important classes of examples are
exhibited. Among these one has: rational normal scrolls, some
Veronese fibrations, some Veronese embeddings of the plane,
defective surfaces, del Pezzo surfaces, etc. \par

With all the above apparatus at hand, the natural question is to
look for classification theorems for $\mathcal{M}^{k}$,
$\mathcal{MA}^{k+1}_{k-1}$, and
$\mathcal{OA}^{k+1}_{k-1}$--varieties. This turns out to be a very
intriguing but considerably difficult question to answer. Indeed
the problem is non trivial even in the case of curves, considered
in \S \ref{curves}: the classification theorem here, which follows
by results of Catalano--Johnson, is that a curve is an
$\mathcal{MA}^{k+1}_{k-1}$ or an
$\mathcal{OA}^{k+1}_{k-1}$--variety if and only if it is a
rational normal curve (see Theorem \ref{curvess}). Our proof is a
slight variation of Catalano--Johnson's argument.  The
classification of $\mathcal{OA}^{2}_{0}$--varieties, also called
$OADP$--varieties, which means {\it varieties with one apparent
double point}, is a classical problem. The case of
$OADP$--surfaces goes back to Severi \cite{Se}, whereas examples
and general considerations concerning the higher dimensional case
can be found in papers by Edge \cite{Edge} and Bronowski
\cite{Br1}. This latter author came to the consideration of this
problem studying extended forms of the Waring problem for
polynomials. Severi's incomplete argument has been recently fixed
by the second author \cite{Ru}, and a different proof can be found
in \cite{CMR}, where one provides the full classification of
$OADP$--threefolds in $\p^7$. Finally an attempt of classification
of $\mathcal{OA}^{k+1}_{k-1}$--surfaces is again due to Bronowski
\cite {Br}, whose approach, based on his aforementioned unproved
conjecture, was certainly not rigorous and led him, by the way, to
an incomplete list.\par

The problem we started from, and which actually was the original
motivation for this paper, was to verify and justify Bronowski's
classification theorem of $\mathcal{OA}^{k+1}_{k-1}$--surfaces,
without, unfortunately, having the possibility of fully relying on
his still unproven conjecture. It was in considering this question
that we understood we had to slightly change our viewpoint and
first look at a different kind of problem. This leads us to the
second theme of the present paper, i.e. linear system on surfaces,
which occupies \S \ref {CastEnr}. We discovered in fact that the
classification of $\mathcal{MA}^{k+1}_{k-1}$ and
$\mathcal{OA}^{k+1}_{k-1}$--surfaces is closely related to a
beautiful classical theorem of Castelnuovo \cite{Ca} and Enriques
\cite{En} (see Theorem \ref{Enr}) which gives an upper bound for
the dimension of a linear system $\cal L$ of curves of given
geometric
 genus
on a surface $X$, and classifies those pairs $(X,{\cal L})$ for
which the bound is attained. Of course Castelnuovo--Enriques'
theorem has to do with the intrinsic birational geometry of
surfaces. However, if one looks at the hyperplane sections linear
systems, it becomes a theorem in projective geometry and our
remark was that Castelnuovo--Enriques' list of extremal cases
consisted of some $k$--defective surfaces and of
$\mathcal{MA}^{k+1}_{k-1}$ and
$\mathcal{OA}^{k+1}_{k-1}$--surfaces for some $k$. It became then
apparent to us that there should have been a relationship between
minimality properties of secant varieties encoded in the
$\mathcal{M}^{k}$, $\mathcal{MA}^{k+1}_{k-1}$, and
$\mathcal{OA}^{k+1}_{k-1}$--properties and the
Castelnuovo--Enriques' maximality conditions on the dimension of
the hyperplane sections linear system. The relation between the
two items was underlined, in our view, by the fact that
Castelnuovo and Enriques' beautiful original approach was based on
iterated applications of tangential projections, a technique that,
as we indicated above, enters all the time in the study of secant
varieties. In fact, we do not reproduce here
Castelnuovo--Enriques' original argument, which, based on the
technical Proposition \ref{tangent}, is however hidden, as we will
explain in a moment, in the proof of our classification theorems
of  $\mathcal{M}^{k}$, $\mathcal{MA}^{k+1}_{k-1}$, and
$\mathcal{OA}^{k+1}_{k-1}$--surfaces given in \S\S \ref{oasurf}
and \ref {msurf}. We preferred instead to give an intrinsic,
birational geometric, proof of Castelnuovo--Enriques' theorem,
which enables us to prove a slightly  more general statement than
the original one and is also useful for extensions, like our
Theorem \ref{3k+2}, in which we classify those smooth surfaces in
projective space such that their hyperplane linear system has
dimension {\it close} to Casteluovo--Enriques' upper bound. The
Castelnuovo-Enriques' upper bound \eqref{grado} for smooth
irreducible curves is essentially the main result of \cite{hart},
see Corollary 2.4, Theorem 3.5 and Theorem 4.1 of {\it loc. cit.},
where the classification of the extremal cases is not considered.
 Our simple and short
proof, which we hope has some independent interest, relies on an
application of  Mori's Cone Theorem, namely Proposition \ref{adj},
which has an independent interest and says that given a pair
$(X,D)$, where $X$ is a smooth, irreducible, projective surface,
and $D$ is a nef divisor on it, one has that $K+D$ is also nef,
unless one of the following facts occurs: either $(X,D)$ is not
{\it minimal}, i.e. there is an exceptional curve of the first
kind $E$ on $X$ such that $D\cdot E=0$, or $(X,D)$ is a $h$--{\it
scroll}, with $h\leq 1$, i.e. there is a rational curve $F$ on $X$
such that $F^2=0$ and $D\cdot F=h$, or $(X,D)$ is a $d$--{\it
Veronese}, with $d\leq 2$, i.e. $X=\p^2$ and $D$ is a curve of
degree $d\leq 2$. A slightly more general version of this last
result, in the case $D$ irreducible (smooth) curve, was obtained
by Iitaka, see \cite{Iitaka}, and revised from the above point of
view of the Cone Theorem by Dicks, see \cite{Dicks} Theorem 3.1.
For weaker results of the same type, concerning the case $D$
ample, see  for example \cite{io2}. It should be stressed that, as
indicated in Castelnuovo's paper \cite{Ca2}, one can push these
ideas further, thus giving suitable upper bounds for the dimension
of certain linear systems on scrolls, or equivalently on the
degree of curves on scrolls as in \cite{hart}, Theorem 2.4 and
Corollary 2.5. This has been done already, in an independent way
also in \cite{Re}, but we hope to return on these matters in the
future since we believe that some of the results in \cite{Ca2},
see also \cite{hart} sections 2 and 3, and in \cite{Re} can be
slightly improved and perhaps related to projective geometry in
the spirit of the present paper.
\par

As we said, in \S \S \ref{oasurf} and \ref {msurf} we come back to
the classification of $\mathcal{MA}^{k+1}_{k-1}$ and
$\mathcal{OA}^{k+1}_{k-1}$--surfaces. Using the machinery of
tangential projections and degeneration of projections we discover
that the surfaces in question are either extremal with respect to
Castelnuovo--Enriques' bound or they are close to be extremal, so
that their classification can be at this point accomplished using
the results of \S \ref {CastEnr}. Finally in \S \ref {gensev} we
prove, using the same ideas, a result, namely Theorem
\ref{Bronowskidifettivo}, which is a wide generalization of the
famous theorem of Severi's saying that the Veronese surface in
$\p^5$ is the only defective surface which is not a cone. \par

In conclusion we would like to mention that, though the above
classification results for $\mathcal{M}^{k}$,
$\mathcal{MA}^{k+1}_{k-1}$, and
$\mathcal{OA}^{k+1}_{k-1}$--varieties are quite satisfactory and
conclusive in low dimensions, i.e. for curves and surfaces, quite
a lot of room is left open for the higher dimensional case, where,
except for the aforementioned result of \cite{CMR}, nothing, to
the best of our knowledge, is known. We hope the ideas presented
in this paper will be useful in this more general context too.
Another interesting direction of research is to try to extend to
higher dimensional varieties Castelnuovo--Enriques' results in \S
\ref{CastEnr}. This question is also widely open. The adjunction
theoretical approach that we use in the surface case can in
principle be extended, but it is not clear whether it leads to
anything really useful. On the other hand Castelnuovo--Enriques
tangential projection approach, in order to work, has to be
modified, since one needs to make projections from osculating,
rather than tangent, spaces. An interesting suggestion in this
direction comes from the beautiful comments of Castelnuovo's to
\cite{Ca} in the volume of collected papers \cite{Cacp}, pp.
186--188. However osculating projections present serious technical
problems which make Castelnuovo's suggestion rather hard to be
pursued. On the other hand the specific problem which Castelnuovo
was considering in his comments in \cite{Cacp}, pp. 186--188, i.e.
the classification of linear systems of rational surfaces in
$\p^3$, has been recently successfully addressed by various
authors, in particular by Mella \cite{Me}, by using Mori's
program. The interplay between intrinsic birational geometry, i.e.
Mori's program, and extrinsic projective geometry, i.e. osculating
projections and relations with secant varieties, is a very
promising, uncharted territory to be explored.\medskip

\noindent {\bf Acknowledgements}: this research started during a
visit of the second author at the University of Roma Tor Vergata
during the months October--December 2002, supported by the
Istituto Nazionale di Alta Matematica "F. Severi". The first
author wants to thank the IMAR of Bucharest for the warm
hospitality in the month of February 2003, when he gave a series
of lectures on the subject of the present paper. The remarks of
the audience, especially of proff. V. Brinzanescu, P. Ionescu and
M. Mendes Lopes, who was also visiting IMAR at the same time, have
been stimulating and precious for him. Both authors finally want
to thank proff. M. Beltrametti and L. B\u adescu for various
discussions and helpful comments about the proof of the results in
\S \ref {CastEnr}. The second author was partially supported by
CNPq, grant nr. 300961/2003-0, and by PRONEX-Algebra Comutativa e
Geometria Algebrica.

\section {Notation and preliminary results}\label {notation}

\subsection {} Let $X\subseteq\p^r$ be a projective scheme over $\C$. We
will denote by $\deg(X)$ the {\it degree} of $X$,  by $\dim(X)$
the {\it dimension} of $X$, by $\codim(X)=r-\dim(X)$ its {\it
codimension} and by $(X)_{\rm red}$ the reduced subscheme
supported by $X$. We will mainly consider the case in which $X$ is
a reduced, irreducible variety. \par

If $Y\subset \p^r$ is a subset, we denote by $<Y>$ the span of
$Y$. We will say that $Y$ is {\it non--degenerate} if
$<Y>=\p^r$.\par

\subsection {} Let $X\subseteq\p^r$ be a reduced, irreducible variety of
dimension $n$. If $x\in X$ we will denote by $C_{X,x}$ the {\it
tangent cone} to $x$ at $X$, which is an $n$--dimensional cone
with vertex at $x$. Note that $C_{X,x}$ has a natural structure of
a subscheme of $\p^r$. We will denote by $\mult_x(X)$ the {\it
multiplicity} of $X$ at $x$. One has $\mult_x(X)=\deg (C_{X,x})$
and $X$ is a {\it cone} if and only if $X$ has some point $x$ such
that $\mult_x(X)=\deg(X)$. In this case $x$ is a {\it vertex} of
$X$ and we will denote by ${\rm Vert}(X)$ the set of vertices of
$X$, which is a linear subspace contained in $X$. It is well known
that:

\begin{equation}\label{cone} {\rm Vert}(X)= \bigcap_{x\in
X}T_{X,x}\end{equation}

If $x$ is a smooth point of $X$, then $C_{X,x}$ is an
$n$--dimensional linear subspace of $\p^r$, i.e. the {\it tangent
space} to $X$ at $x$, which we will denote by $T_{X,x}$.\par

\subsection {}\label{subsecsec} Let $k$ be a non--negative integer and
let $S^{k}(X)$ be the $k$--{\it secant variety} of $X$, i.e. the
Zariski closure in $\p^r$ of the set:\medskip

\centerline {$\{ x\in \p^r: x$  lies in the span of $k+1$
independent points of $X\}$}\medskip

Of course $S^0(X)=X$, $S^r(X)=\p^r$ and $S^{k}(X)$ is empty if
$k\geq r+1$. We will write $S(X)$ instead of $S^1(X)$ and we will
assume $k\leq r$ from now on. \par

Let ${\rm Sym}^{h}(X)$ be the $h$--th symmetric product of $X$.
One can consider the {\it abstract $k$--th secant variety} $S^k_X$
of $X$, i.e. $S^{k}_X\subseteq {\rm Sym}^k(X)\times\p^r$ is the
Zariski closure of the set of all pairs $([p_0,...,p_k],x)$ such
that $p_0,...,p_k\in X$ are linearly independent points and $x\in
<p_0,...,p_k>$. One has the surjective map $p_X^k:S^k_X\to
S^k(X)\subseteq\p^r$, i.e. the projection to the second factor.
Hence:\medskip

\begin{equation} \label{defect} s^{(k)}(X):= \dim (S^{k}(X))\leq
\min\{r,\dim(S^k_X)\}= \min\{r,n(k+1)+k\}\end{equation}

We will denote by $h^{(k)}(X)$ the codimension of $S^{k}(X)$ in
$\p^r$, i.e. $h^{(k)}(X):=r- s^{(k)}(X)$.\par

The right hand side of (\ref {defect}) is called the {\it expected
dimension} of $S^{k}(X)$ and will be denoted by $\sigma^{(k)}(X)$.
One says that $X$ has a $k$-{\it defect}, or is $k$-{\it
defective}, or is {\it defective of index} $k$ when strict
inequality holds in (\ref {defect}). One says that:

$$\delta_k(X):=\sigma^{(k)}(X)-s^{(k)}(X)$$

\noindent is the $k$--{\it defect} of $X$.\par

Notice that the general fibre of $p_X^k$ is pure of dimension
$(k+1)n+k-s^{(k)}(X)$, which equals $\delta_k(X)$ when $r\geq
n(k+1)+k$. We will denote by $\mu_k(X)$ the number of irreducible
components of this fibre. In particular, if $s^{(k)}(X)=(k+1)n+k$,
then $p_X^k$ is generically finite and $\mu_k(X)$ is the degree of
$p_X^k$, i.e. it is the number of $(k+1)$--secant $\p^k$'s to $X$
passing through the general point of $S^k(X)$.\par

If $s^{(k)}(X)=(k+1)n+k$, we will denote by $\nu_k(X)$ the number
of $(k+1)$--secant $\p^k$'s to $X$ meeting the general
$\p^{h^{(k)}(X)}$ in $\p^r$. Of course one has:

\begin{equation} \label {munu} \nu_k(X)=\mu_k(X)\cdot
\deg(S^k(X))\end{equation}

\noindent and therefore:

\begin{equation} \label {mu=nu} \nu_k(X)=\mu_k(X) \quad {\rm if}\quad
r=s^{(k)}(X)=(k+1)n+k. \end{equation}

\subsection {} Let $X\subset\p^r$ be an
irreducible, projective variety. Let $k$ be a positive integer and
let $p_1,...,p_k$ be general points of $X$. We denote by
$T_{X,p_1,...,p_k}$ the span of $T_{X,p_i}, i=1,...,k$.\par

If $X\subset\p^r$ is a projective variety, Terracini's Lemma
describes the tangent space to $S^k(X)$ at a general point of it
(see \cite {Terr1} or, for modern versions, \cite {Adl}, \cite
{CC}, \cite {Dale1}, \cite {Zak}):

\begin{Theorem}\label{terracini} (Terracini's Lemma) Let $X\subset\p^r$
be an irreducible, projective variety. If $p_0,...,p_k\in X$ are
general points and $x\in <p_0,...,p_k>$ is a general point, then:

$$T_{S^k(X),x}=T_{X,p_0,...,p_k}.$$

If $X$ is $k$--defective, then the general hyperplane $H$
containing $T_{X,p_0,...,p_k}$ is tangent to $X$ along a variety
$\Sigma_{p_0,...,p_k}$ of pure, positive dimension $n_k(X)$
containing $p_0,...,p_k$. Moreover one has:

$$k\leq \dim(<\Sigma_{p_0,...,p_k}>)\leq kn_k(X)+k+n_k(X)-\delta_k(X).$$

\end{Theorem}

Consider the projection of $X$ with centre $T_{X,p_1,...,p_k}$. We
call this a {\it general $k$--tangential projection} of $X$, and
we will denote it by $\tau_{X,p_1,...,p_k}$ or simply by
$\tau_{X,k}$. We will denote by $X_k$ its image. By Terracini's
lemma, the map $\tau_{X,k}$ is generically finite to its image if
and only if $s^{(k)}(X)=(k+1)n+k$. In this case we will denote by
$d_{X,k}$ its degree.\par

In the same situation, the projection of $X$ with centre the space
$<p_1,...,p_k>$ is called a {\it general $k$--internal projection}
of $X$, and we will denote it by $t_{X,p_1,...,p_k}$ or simply by
$t_{X,k}$. We denote by $X^k$ its image. We set $X_0=X^0=X$.
Notice that the maps $t_{X,k}$ are birational to their images as
soon as $k<r-n=\codim(X)$.  \par

Sometimes we will use the symbols $X_k$ [resp. $X^k$] for
$k$--tangential projections [resp. $k$--internal projections]
relative to {\it specific}, rather than {\it general}, points. In
this case we will explicitly specify this, thus we hope no
confusion will arise for this reason.

\subsection{}\label{subsecwd} We recall from \cite {CC} the definition of
a {\it $k$--weakly defective} variety, i.e. a variety
$X\subset\p^r$ such that if $p_0,...,p_k\in X$ are general points,
then  the general hyperplane $H$ containing $T_{X,p_0,...,p_k}$ is
tangent to $X$ along a variety $\Sigma_{p_0,...,p_k}$ of pure,
positive dimension $n_k(X)$ containing $p_0,...,p_k$. By
Terracini's lemma, a $k$--defective variety is also $k$--weakly
defective, but the converse does not hold in general (see \cite
{CC}).

\begin{Remark} \label{0wd} {\rm A curve is never
$k$--weakly defective for any $k$. A variety is $0$--weakly
defective if and only if its dual variety is not a hypersurface.
In the surface case this happens if and only if the surface is
{\it developable}, i.e. if and only if the surface is either a
cone or the tangent developable to a curve.}\end{Remark}

The two next results are consequences of Theorem 1.4 of \cite{CC}
that we partially recall here.

\begin{Theorem} Let $X\subset\p^r$ be an irreducible, projective,
non--degenerate variety of dimension $n$. Assume $X$ is not
$k$--weakly defective for a given $k$ such that $r\geq
(n+1)(k+1)$. Then, given $p_0,\ldots, p_{k}$ general points on
$X$, the general hyperplane $H$ containing $T_{X,p_0,...,p_k}$ is
tangent to $X$ only at $p_0,\ldots, p_{k}$. Moreover such a
hyperplane $H$ cuts on $X$ a divisor with ordinary double points
at $p_0,\ldots, p_{k}$.\label{CC}
\end{Theorem}

The first consequence we are interested in is the following:

\begin{Lemma} Let $X\subset \p^r$ be an irreducible, projective,
non--degenerate variety of dimension $n$, which is not $k$--weakly
defective for a fixed $k\geq 1$ such that $r\geq (k+1)(n+1)$. Then
a general $k$--tangential projection of $X$ is birational to its
image, i.e. $d_{X,k}=1$. In particular, if $r\geq 2n+2$, the
general tangential projection of $X$ is birational to its
image.\label{birational}
\end{Lemma}

\begin{proof} Since $X$ is not $k$-weakly defective, it is not
$l$--defective for all $l\leq k$. Thus we have
$s^{(l)}(X)=(l+1)n+l$ for all $l\leq k$, so that by Terracini's
lemma $\tau_{X,p_1,...,p_l}$ is generically finite onto $X_l$ for
every $l\leq k$ and $p_1,...,p_l$ general points on $X$. In
particular this is true for $l=k$. \par

Suppose now that $d_{X,k}>1$. Then, given a general point $p_0\in
X$ there is a point $q\in X\setminus (T_{X,p_1,...,p_{k}}\cap X)$,
$q\neq p_0$, such that
$\tau_{X,p_1,...,p_{k}}(p_0)=\tau_{X,p_1,...,p_{k}}(q):=x\in
X_{k}$. This would imply that $T_{X,p_0,p_1,...,p_k}$ and
$T_{X,q,p_1,...,p_{k}}$ coincide, since both these spaces project
via $\tau_{X,p_1,...,p_{k}}$ onto $T_{X_k,x}$. In particular, the
general hyperplane tangent to $X$ at $p_0,p_1,\ldots,p_k$ is also
tangent at $q$. This contradicts Theorem \ref{CC}.
\end{proof}

We also note that Terracini's lemma and Theorem \ref{CC} imply
that:

\begin{Proposition}\label{wd} Let $X\subset\p^r$ be an
irreducible, projective variety which is not $k$--weakly
defective. If $r\geq (n+1)(k+1)$, then
$\mu_k(X)=1$.\end{Proposition}

In the sequel we will also need the following technical:

\begin{Proposition} \label{tangent} Let $X\subset \p^r$ be a smooth,
irreducible, projective, non--degenerate surface, which is not
$(k-1)$-weakly defective for a fixed $k\geq 1$ such that $r\geq
3k+2$. Let $p_1, \ldots, p_{k}\in X$ be general points and assume
that the linear system $\mathcal L$ of hyperplane sections of $X$
tangent at $p_1, \ldots, p_{k}$ has a not empty fixed part
$F=\sum_{i=1}^h n_i\Gamma_i$, with $\Gamma_i$ distinct,
irreducible curves and $n_i>0$, for all $i=1,...,h$. Let $\mathcal
M$ be the movable part of $\mathcal L$ and let $M$ be its general
curve. Then $F$ is reduced, i.e. $n_i=1$ for all $i=1,...,h$ and:

\begin{itemize}

\item[(i)] either $h=1$, $F$ is a smooth, rational curve
containing $p_1, \ldots, p_{k}$, whereas $\mathcal M$ has simple
base points at $p_1, \ldots, p_{k}$ and $M\cdot F=k$, hence $M\in
\mathcal M$ general meets $F$ transversally at $p_1, \ldots,
p_{k}$ and nowhere else;\par

\item[(ii)] or $h=k$, $\Gamma_i$ is a smooth, rational curve
containing $p_i$ for $i=1,...,k$, $\Gamma_i\cap\Gamma_j=\emptyset$
if $1\leq i< j\leq k$, $\mathcal M$ has simple base points at
$p_1, \ldots, p_{k}$ and $M\cdot \Gamma_i=1$, hence $M\in \mathcal
M$ general meets $\Gamma_i$ transversally at  $p_i$ and nowhere
else, for all $i=1,...,k$. \end{itemize}

Moreover, if $r\geq 3k+3$ and if the general $k$--tangential
projection $X_{k}$ of $X$, has rational hyperplane sections, then
the general curve $M\in \mathcal M$ is rational. \end{Proposition}

\begin{proof} Let $C$ be a general curve in $\mathcal L$, so that $C=F+M$.
By Theorem \ref {CC}, we know that $C$ has nodes at $p_1, \ldots,
p_{k}$ and is otherwise smooth. This implies that:

\begin{itemize}

\item $F$ is reduced;

\item all the curves $\Gamma_i$, $i=1,...,h$, are smooth off $p_1,
\ldots, p_{k}$, where they can have at most nodes;

\item $\Gamma_i$ and $\Gamma_j$, for $1\leq i< j\leq h$, may
intersect only at some of the points $p_1, \ldots, p_{k}$, where
only two of them may meet transversally;

\item $M$ is smooth off $p_1, \ldots, p_{k}$ where it can have at
most nodes,  and may intersect the curves $\Gamma_i$ only at $p_1,
\ldots, p_{k}$, where it may meet only one of them transversally;

\item if the point $p_i$, $i=1,...,k$, is a node for a curve
$\Gamma_j$, $i=1,...,h$, then it does not belong neither to $M$,
nor to $\Gamma_j$, $j\not=i$;\par

\item if the point $p_i$, $i=1,...,k$, is a node for $M$, then it
does not belong to $F$;\par

\item if the point $p_i$, $i=1,...,k$, is a smooth point for a
curve $\Gamma_j$, $i=1,...,h$, then it  belongs either to $M$, or
to a curve $\Gamma_j$, $j\not=i$, but not to both.

\end{itemize}

We prove the assertion in various steps.\medskip

\begin{Claim} \label{claima} Every irreducible component $\Gamma_i$ of $F$
contains some of the points  $p_1, \ldots, p_{k}$.
\end{Claim}\medskip

Otherwise we would have $\Gamma_i\cap \overline
{C-\Gamma_i}=\emptyset$, and $C$ would be disconnected, a
contradiction since it is very ample on $X$.\medskip

\begin{Claim} \label{claimb} $F$ contains all the points $p_1, \ldots,
p_{k}$.\end{Claim}\medskip

In fact, if $p_1\not\in F$, then, by changing the role of the
points $p_1, \ldots, p_{k}$, none of the points $p_1, \ldots,
p_{k}$ is in $F$, contradicting Claim \ref {claima}.\medskip

\begin{Claim}\label{claimd} $F$ is smooth.\end{Claim}\medskip

We know $F$ can be singular only at some of the points
$p_1,...,p_{k}$. Suppose this is the case. Then by symmetry, it is
singular at any one of the points in question. But then we would
have $M\cap F=\emptyset$, which leads to a contradiction as above.
\medskip

\begin{Claim} \label{claimc} Let $\Gamma_1$ be the irreducible component
of $F$ through $p_1$. Then either also $p_2,...,p_{k}\in
\Gamma_1$, or none of the points $p_2,...,p_{k}$ lies on
$\Gamma_1$. In the former case $\Gamma_1=F$. In the latter each of
the points $p_i, i=1,...,k$, belongs to one and only one component
$\Gamma_i$ of $F$. \end{Claim}\medskip

Suppose $\Gamma_1$ contains $p_1,...,p_i$, with $1<i<k$. By
changing the role of the points $p_1, \ldots, p_{k}$, any $i$
among the points  $p_1,...,p_{k}$ lie on some irreducible
component of $F$. Then $F$ would be singular, contradicting Claim
\ref {claimd}. This proves the first part of the Claim.
\par

Assume $p_1,...,p_{k}\in \Gamma_1$. Then Claim \ref {claima} and
Claim \ref {claimd} imply that $F=\Gamma$. Suppose instead only
$p_1$ lies on $\Gamma$. Then by changing the role of the points
$p_1, \ldots, p_{k}$, each of the other points $p_i, i=2,...,k$,
also lies on one and only one component of $F$.\medskip

\begin{Claim}\label{claime} Every irreducible component $\Gamma_i$ of $F$
is rational.\end{Claim}

By projecting $X$ from  $T_{X,p_1, \ldots, p_{k-1}}$, we get an
irreducible surface $X_{k-1}\subset \p^{r-3k+3}$, with $r-3k+3\geq
5$, which is birational to $X$ by Lemma \ref{birational} and which
is not $0$--weakly defective. Let $q$ be the image on $X_{k-1}$ of
a general point $p_k$ of $X$. Notice that the general tangent
hyperplane section to $X_{k-1}$ at $q$, which is the image of $C$,
is reducible containing $M'$, the image of $M$, and $\Gamma'$, the
image of $\Gamma_k$, both passing through $q$. Notice that $M'$ is
the movable part of the linear system of hyperplane sections of
$X_{k-1}$ tangent at $q$, whereas $\Gamma'$ is the fixed part.
Then $X_{k-1}$ is either the Veronese surface in $\p^5$ or a
non--developable scroll over a curve (see for instance
\cite{mez}). Hence $\Gamma'$ is rational. Since
$\tau_{X,p_1,...,p_{k-1}}$ is birational by Lemma
\ref{birational}, then $\Gamma_k$ is birational to $\Gamma'$, and
is therefore rational. If $\Gamma_k=F$ there is nothing else to
prove. Otherwise, by changing the role of the points $p_i$, we see
that $\Gamma_i$ is rational for any $i=1,...,k$.\medskip

The above claims imply (i) and (ii). As for the last assertion, it
 follows from Lemma \ref{birational}. \end{proof}

\subsection{}  If $X, Y\subset \p^r$ are closed subvarieties we denote by
$J(X,Y)$ the {\it join} of $X$ and $Y$, i.e. the Zariski closure
of the union of all lines $<x,y>$, with $x\in X, y\in Y, x\not=
y$. If $X$ is a linear subspace, then $J(X,Y)$ is the cone over
$Y$ with vertex $X$. With this notation, for every $k\geq 1$ one
has:

\begin{equation} \label{join}
S^k(X)=J(S^l(X),S^h(X))\end{equation}

\noindent if $l+h=k-1$, $l\geq 0$, $h\geq 0$. \par

We record the following:

\begin{Lemma}\label{projoin} Let $X, Y\subset \p^r$ be closed,
irreducible, subvarieties and let $\Pi$ be a linear subspace of
dimension $n$ which does not contain either $X$ or $Y$. Let $\pi:
\p^r\map \p^{r-n-1}$ be the projection from $\Pi$ and let $X', Y'$
be the images of $X, Y$ via $\pi$. Then:

$$\pi(J(X,Y))=J(X',Y').$$

\noindent In particular, if $\Pi$ does not contain $X$, then for
any non--negative inter $k$ one has:

$$\pi(S^k(X))=S^k(X').$$\end{Lemma}

\begin{proof} It is clear that $\pi(J(X,Y))\subseteq J(X',Y')$. Let $x'\in
X', y'\in Y'$ be general points. Then there are $x\in X, y\in Y$
such that $\pi(x)=x', \pi(y)=y'$. Thus $\pi(<x,y>)=<x',y'>$,
proving that $J(X',Y')\subseteq \pi(J(X,Y))$, i.e. the first
assertion. The rest of the statement follows by (\ref {join}) with
$l=0$, by making induction on $k$.\end{proof}

The following lemma is an application of Terracini's lemma:

\begin{Lemma}\label {kl} Let $X\subset\p^r$ be an
irreducible, projective variety. For all $i=1,...,k$ one has:
$$h^{(k-i)}(X_i)=h^{(k)}(X)$$
whereas for all $i\geq 1$ one has:
$$h^{(k)}(X^i)=\max\{0,h^{(k)}(X)-i\}.$$\end{Lemma}

\begin{proof} Let $p_0,...,p_k\in X$ be general points. Terracini's lemma
says that $T_{X,p_0,...,p_k}$ is a general tangent space to
$S^k(X)$ and that its projection from $T_{X,p_{k-i+1},...,p_k}$ is
the general tangent space to $S^{k-i}(X_i)$. This implies the
first assertion.\par

To prove the second assertion, note that it suffices to prove it
for $i< h^{(k)}(X)$. Indeed, if $i\geq h^{(k)}(X)$ then, by Lemma
\ref {projoin} one has $h^{(k)}(X^i)=0$ since already
$h^{(k)}(X^{h^{(k)}})=0$. Thus, suppose $i< h^{(k)}(X)$.  Let
$p_0,...,p_k\in X$ be general points and take $i$ general points
$q_1,...,q_i$ in $X\setminus(X\cap T_{X,p_0,...,p_k})$. Then the
projection of $T_{X,p_0,...,p_k}$ from $<q_0,...,q_i>$ is the
tangent space to $S^k(X^i)$. Furthermore $i< h^{(k)}(X)$ yields
$<q_0,...,q_i>\cap T_{X,p_0,...,p_k}=\emptyset$. This implies the
second assertion. \end{proof}

\subsection {}\label{RNS} Let $0\leq a_1\leq a_1\leq ...\leq a_n$ be
integers and set $\p(a_1,...,a_n):=\p({\cal
O}_{\p^1}(a_1)\oplus...\oplus{\cal O}_{\p^1}(a_n))$.  We will
denote by $H$ a divisor in $|{\cal O}_{\p(a_1,...,a_n)}(1)|$ and
by $F$ a fibre of the structure morphism $\pi: \p(a_1,...,a_n)\to
\p^1$. Notice that the corresponding divisor classes, which we
still denote by $H$ and $F$, freely generate ${\rm
Pic}(\p(a_1,...,a_n))$.\par

Set $r=a_1+...+a_n+n-1$ and consider the morphism:

$$\phi:=\phi_{|H|}: \p(a_1,...,a_n)\to \p^r$$

\noindent whose image we denote by $S(a_1,...,a_n)$. As soon as
$a_n>0$, the morphism $\phi$ is birational to its image. Then the
dimension of $S(a_1,...,a_n)$ is $n$ and its degree is
$a_1+...+a_n=r-n+1$, thus $S(a_1,...,a_n)$ is a {\it rational
normal scroll}, which is smooth if and only if $a_1>0$. Otherwise,
if $0=a_1=...=a_i<a_{i+1}$, then $S(a_1,...,a_n)$ is  the cone
over $S(a_{i+1},...,a_n)$ with vertex a $\p^{i-1}$. One uses the
simplified notation $S(a_1^{h_1},...,a_m^{h_m})$ if $a_i$ is
repeated $h_i$ times, $i=1,...,m$.

We will sometimes use the notation $H$ and $F$ to denote the Weil
divisors in $S(a_1,...,a_n)$ corresponding to the ones on
$\p(a_1,...,a_n)$. Of course this is harmless if $a_1>0$, since
then $\p(a_1,...,a_n)\simeq S(a_1,...,a_n)$.
\par

Recall that rational normal scrolls, the Veronese surface in
$\Proj^5$ and the cones on it, and the quadrics,  can be
characterized as those non--degenerate, irreducible varieties
$X\subset \p^r$ in a projective space having minimal degree
$\deg(X)=\codim(X)+1$ (see \cite{EH}).\par

Let $X=S(a_1,\ldots,a_n)\subset\p^r$ be as above. We leave to the
reader to see that:

\begin{equation} \label{intproj} X^1=S(b_1,...,b_n)\quad {\rm
where}\quad \{b_1,...,b_n\}=\{a_1,...,a_n-1\}.\end{equation}

One can also consider the projection $X'$ of $X$ from a general
$\p^{n-1}$ of the ruling of $X$. This is not birational to its
image if $a_1=0$ and one sees that  if  $a_1=...=a_i=0<a_{i+1}$,
then:

\begin{equation} \label{projsub} X'=S(c_1,...,c_{n-i})\quad {\rm
where}\quad
\{c_1,...,c_{n-i}\}=\{a_{i+1}-1,...,a_n-1\}.\end{equation}

A general tangential projection of $X=S(a_1,\ldots,a_n)$ is the
composition of the projection of $X$ from a general $\p^{n-1}$ of
the ruling of $X$ and of a general internal projection of $X'$.
Therefore, by putting (\ref{intproj}) and (\ref{projsub})
together, one deduces that  if  $a_1=...=a_i=0<a_{i+1}$, then:

\begin{equation} \label{tanproj} X_1=S(d_1,...,d_{n-i})\quad {\rm
where}\quad
\{d_1,...,d_{n-i}\}=\{a_{i+1}-1,...,a_n-2\}.\end{equation}

As a consequence we have:

\begin{Proposition}\label{defscroll} Let
$X=S(a_1,\ldots,a_n)\subset\p^r$ be a rational normal scroll as
above. Then:

$$\dim(S^k(X))=\min\{r, r+k+1-\sum_{1\leq j\leq n;\; k\leq a_j}
(a_j-k)\}.$$

In particular, if $r\geq (k+1)n+k$, then $s^{(k)}(X)=(k+1)n+k$ if
and only if $a_1\geq k$.

\end{Proposition}

\begin{proof} It follows by induction using (\ref {tanproj}) and
Terracini's lemma. We leave the details to the reader.\end{proof}

A different proof of the same result can be obtained by writing
the equations of $S^k(X)$ (see \cite{Room} and \cite{CJ2} for this
point of view).

\subsection {} Given positive integers $0<m_1\leq ...\le m_h$ we will
denote by ${\rm Seg}(\p^{m_1},...,\p^{m_h})$, or simply by ${\rm
Seg}(m_1,...,m_h)$ the {\it Segre variety} of type
$(m_1,...,m_h)$, i.e. the image of  $\p^{m_1}\times ...\times
\p^{m_h}$ in $\Proj^r$, $r=(m_1+1)\cdot\cdot\cdot(m_h+1)-1$, under
the  {\it Segre embedding}. Notice that, if $\p^{m_i}=\p(V_i)$,
where $V_i$ is a complex vector space of dimension $m_i+1$,
$i=1,...,h$, then $\p^r=\p(V_1\otimes ...\otimes V_h)$ and ${\rm
Seg}(m_1,...,m_h)$ is the set of equivalence classes of
indecomposable tensors in $\p^r$. We use the shorter notation
${\rm Seg}(m_1^{k_1},...,m_s^{k_s})$ if $m_i$ is repeated $k_i$
times, $i=1,...,s$.\par

 Recall that ${\rm Pic}(\p^{m_1}\times ...\times \p^{m_h})\simeq {\rm
Pic}({\rm Seg}(m_1,...,m_h))\simeq \Z^h$, is freely generated by
the line bundles $\xi_i=pr_i^*({\cal O}_{\Proj^{m_i}}(1))$,
$i=1,...,h$, where $pr_i: \p^{m_1}\times ...\times \p^{m_h}\to
\p^{m_i}$ is the projection to the $i$--th factor. A divisor $D$
on ${\rm Seg}(m_1,...,m_h)$ is said to be of ${\it type}$ $
(\ell_1,...,\ell_h)$ if ${\cal O}_{{\rm
Seg}(m_1,...,m_h)}(D)\simeq \xi_1^{\ell_1}\otimes ... \otimes
\xi_h^{\ell_h}$. The line bundle $\xi_1^{\ell_1}\otimes ...
\otimes \xi_h^{\ell_h}$ on $\p^{m_1}\times ...\times \p^{m_h}$ is
also denoted by ${\cal O}_{\p^{m_1}\times ...\times
\p^{m_h}}(\ell_1,...,\ell_h)$. The hyperplane divisor of ${\rm
Seg}(m_1,...,m_h)$ is of type $(1,...,1)$.\par

It is useful to recall what are the defects of the Segre varieties
${\rm Seg}(m_1,m_2)$ with $m_1\leq m_2$. As above, let $V_i$ be
complex vector spaces of dimension $m_i+1$, $i=1,2$. We can
interpret the points of  $\p(V_1\otimes V_2)$ as the equivalence
classes of all $(m_1+1)\times (m_2+1)$ complex matrices and ${\rm
Seg}(m_1,m_2)={\rm Seg}(\p(V_1),\p(V_2)))$ as the subscheme of
$\p(V_1\otimes V_2)$ formed by the equivalence classes of all
matrices of rank $1$. Similarly $S^{k}({\rm Seg}(m_1,m_2))$ can be
interpreted as the subscheme of  $\p(V_1\otimes V_2)$ formed by
the equivalence classes of all matrices of rank less than or equal
to $k+1$. Therefore $S^k({\rm Seg}(m_1,m_2))=\p(V_1\otimes V_2)$
if and only if $k\geq m_1$. In the case $k<m_1$ one has instead:

$$\codim(S^k({\rm
Seg}(m_1,m_2))=(m_1-k)(m_2-k)$$

\noindent (see  \cite{acgh}, pg. 67). As a consequence one has:

 $$\delta_k({\rm
Seg}(m_1,m_2))=k(k+1)$$

\noindent if $k<m_1\leq m_2$.\par

The degree of $S^k({\rm Seg}(m_1,m_2))$, with $k<m_1\leq m_2$, are
computed by a well known formula by Giambelli \cite{Giambelli},
apparently already known to C. Segre (see \cite{Room}, pg. 42, and
\cite {Fu}, 14.4.9, for a modern reference). The case $k=m_1-1$,
which is the only one we will use later, is not difficult to
compute (see \cite{Harris}, pg. 243) and reads:

$$\deg(S^{m_1-1}({\rm
Seg}(m_1,m_2))=\binom {m_2+1}{m_1}.$$

\subsection{} We will recall now some definition and result due to G.
Kempf, \cite{Kempf}, which we are going to use later.\par

Let $V_1, V_2, V_3$ finite dimensional complex vector spaces. A
pairing:\par

$$\phi:V_1\otimes V_2\to V_3$$

\noindent is said to be {\it 1--generic} if $0\neq v\in V_1$ and
$0\neq u\in V_2$ implies $\phi(v\otimes u)\neq 0$. From a
projective geometric point of view, $\phi$ determines a projection
$\varphi: \p(V_1\otimes V_2)\map \p(V_3)$ and the $1$--genericity
condition translates into the fact that the centre of the
projection $\varphi$ does not intersect ${\rm
Seg}(\p(V_1),\p(V_2))$. \par

If $\phi$ is surjective, then we may regard $\phi$ as specifying a
linear space of linear transformations:

$$V_3^*\subseteq {\rm Hom}(V_1, V_2^*)\simeq V_1^*\otimes V_2^*.$$

\noindent One says that $V_3^*$ is {\it 1--generic} if $\phi$
is.\par

Let $m_i+1=\dim(V_i)$ and suppose $m_1\leq m_2$. For each $k$ such
that $0\leq k\leq m_1$, let $(V_3^*)_k$ be the subscheme of $V_3^*
$ of all matrices in $V_3^*$ with rank less than or equal to
$k+1$, i.e. the scheme--theoretic intersection of $V_3^*$ with the
scheme ${\rm Hom}(V_1, V_2^*)_k$ of all matrices with rank less
than or equal to $k+1$ in ${\rm Hom}(V_1, V_2^*)$. Of course
$(V_3^*)_k$ is a cone, hence it gives rise to a closed subscheme
$\p((V_3^*)_k)$ of $\p(V_3^*)$ which is the scheme theoretic
intersection of $\p(V_3^*)$ with $S^k({\rm
Seg}(\p(V_1^*),\p(V_2^*))$. Notice that the expected codimension
of $\p((V_3^*)_k)$ in $\p(V_3^*)$ is:

$$m_1m_2-k(m_1+m_2)+k^2=\dim(\p(V_1^*\otimes V_2^*))-s^{(k)}({\rm
Seg}(\p(V_1^*),\p(V_2^*))).$$

\noindent This is also the expected codimension of $(V_3^*)_k$ in
$V_3^*$. We can now state Kempf's theorem:

\begin{Theorem}\label{Kempf} If $V_3^*\subseteq
V_1^*\otimes V_2^*$ is $1$--generic, then $(V_3^*)_{m_1-1}$ is
reduced, irreducible and of the expected codimension $m_2-m_1+1$
in $V_3^*$. The same is true for $\p((V_3^*)_{m_1-1})$, whose
degree is $\binom {m_2+1} {m_1}$.
\end{Theorem}

\subsection{} Given positive integers $n,d$, we will denote by $V_{n,d}$
the image of $\p^n$ under the {\it $d$--Veronese embedding} of
$\p^n$ in $\p^{{{n+d}\choose d}-1}$.

\subsection{} \label{gen} If $X$ is a variety of dimension $n$ and $Y$ a
subvariety of $X$, we will denote by ${\rm Bl}_Y(X)$ the blow--up
of $X$ along $Y$. If $Y$ is a finite set $\{x_1,...,,x_n\}$ we
denote the blow--up by ${\rm Bl}_{x_1,...,x_n}(X)$. \par

With the symbol $\equiv$ we will denote the linear equivalence of
divisors on $X$. The symbol $\sim$ will instead denote numerical
equivalence. If ${\mathcal L}$ is a linear systems of divisors on
$X$, of dimension $r$, we will denote by $\phi_{\mathcal L}: X\map
\p^r$ the rational map defined by $\mathcal L$.\par

If $D$ is a divisor on the variety $X$, we denote by $|D|$ the
complete linear series of $D$. If $X\subset \p^r$ is an
irreducible, projective variety, and $D$ is a hyperplane section
of $X$, one says that $X$ is {\it linearly normal} if the linear
series cut out on $X$ by the hyperplanes of $\p^r$ is complete,
i.e. if the natural map:

$$H^0(\p^r,{\mathcal O}_{\p^r}(1))\to H^0(X,{\mathcal O}_{X}(D))$$

\noindent is surjective.\par

If $D$ [resp. ${\mathcal D}$] is a divisor [resp. a line bundle]
on $X$, we will say that $D$ [resp. ${\mathcal D}$] is {\it
effective} if $h^0(X,{\mathcal O}_X(D))>0$ [resp. $h^0(X,{\mathcal
L})>0$]. We will say that $D$ [resp. ${\mathcal D}$] is {\it nef}
if for every curve $C$ on $X$, one has $D\cdot C\geq 0$ [resp.
${\mathcal D}\cdot C\geq 0$]. A nef divisor $D$ [resp. a nef line
bundle ${\mathcal D}$] is {\it big} if $D^n>0$ [resp. ${\mathcal
D}^n>0$].

\subsection {} \label{gen'} Let $X$ be a smooth, irreducible, projective
surface. As customary, we will use the following notation
$q:=q(X):=h^1(X,\O_X)$ for the {\it irregularity},
$\kappa:=\kappa(X)$ for the {\it Kodaira dimension} of $X$. We
will denote by $K:=K_X$ a {\it canonical divisor} on $X$ and, as
usual, $p_g:=p_g(X):=h^0(X,{\mathcal O}_X(K))$ is the {\it
geometric genus}.
\par

If $C$ is a curve on $X$, it will be called a $(-n)$-{\it curve},
if $C\simeq \Proj^1$ and $C^2=-n$. Recall that a famous theorem of
Castelnuovo's identifies the $(-1)$--curves as the exceptional
divisors of blow--ups.\par

Let $D$ be a Cartier divisor on an irreducible, projective surface
$X$. We denote by $p_a(D)$ the arithmetic genus of $D$. We will
say that $D$ is a {\it curve} on $X$ if it is effective. If $D$ is
reduced curve on $X$, we will consider $p_g(D)$ the {\it geometric
genus} of $D$, i.e. the arithmetic genus of the normalization of
$D$. \par

A curve $D$ on $X$ will be called $m$--connected if for every
decomposition $D=A+B$, with $A,B$ non--zero curves on $X$, one has
$A\cdot B\geq m$. If $D$ is $1$--connected one has $h^0(D,\mathcal
O_D)=1$ and $h^1(D,\mathcal O_D)=p_a(D)\geq 0$ (see \cite {BPV}).
If $D$ is a big and nef curve on $X$, then $D$ is $1$--connected
(see \cite {ML}, Lemma (2.6)).\par

If $X$ is smooth, we will say that the pair $(X,D)$ is:

\begin{itemize}

\item {\it effective}  [resp. {\it nef}, {\it big}, {\it ample},
{\it very ample}] if $D$ is such;

\item {\it minimal} if there is no $(-1)$-curve $C$ on $X$ such
that $D\cdot C=0$;

\item a $h$--{\it scroll}, with $h\geq 0$ an integer, if there is
a smooth rational curve $F$ on $X$ such that $F^2=0$ and $D\cdot
F=h$;

\item  a {\it del Pezzo pair} if $K\sim -D$ and $(X,D)$ is big and
nef.
\end{itemize}

A $1$--scroll will be simply called a {\it scroll}.

Notice that if $(X,D)$ is a del Pezzo pair, then $X$ is rational
and $K\equiv -D$. Indeed $-K$ is nef and big, thus
$\kappa(X)=-\infty$ and $q=h^1(X,\mathcal O_X)=h^1(X,\mathcal
O_X(K-K))=0$ by Ramanujam's vanishing theorem (see \cite
{Ram}).\par

If $\mathcal L$ is a linear system on $X$ and $D\in \mathcal L$ is
its general divisor, we will say that $(X,\mathcal L)$ is {\it
nef, big, ample, minimal, a $h$--scroll} etc. if $(X,D)$ is such.
One says that $(X,\mathcal L)$ is {\it very ample} if
$\phi_\mathcal L$ is an isomorphism of $X$ to its image. \par

Suppose the linear system $\mathcal L$ has no fixed curve and the
general curve in $\mathcal L$ is irreducible. Then, by blowing up
the base points of $\mathcal L$, we see that there is a unique
pair $(X',\mathcal L')$, where $X'$ is a surface with a birational
morphism $f: X'\to X$ and a $\mathcal L'$ is linear system on $X'$
such that:

\begin{itemize}

\item $\mathcal L'$ is the strict transform of $\mathcal L$ on
$X'$;

\item $\mathcal L'$ is base point free, and therefore its general
curve $D'$ is smooth and irreducible;

\item $\mathcal L'$ is $f$--{\it relatively minimal}, i.e. if $E$
is a $(-1)$--curve on $X'$ such that $D'\cdot E=0$ then $E$ is not
contracted by $f$.

\end{itemize}

We will call the pair $(X',\mathcal L')$ the {\it resolution} of
the pair $(X,D)$.

If $X\subseteq \p^r$ is an irreducible, projective surface, one
considers $f: X'\to X\subseteq \p^r$ a minimal desingularization
of $X$ and $\mathcal L$ the linear system on $X'$ such that
$f=\phi_{\mathcal L}$. The pair $(X',\mathcal L)$ is big, nef and
minimal. One says that $X$ is a {\it scroll} if the pair
$(X',\mathcal L)$ is a scroll.\par

If $X\simeq \p^2$ and $R$ is a line, the pair $(X,D)$ with
$D\equiv dR$ will be called a $d$--{\it Veronese pair}. If
$X=\FF_a:=\p(0,a)$ is the Hirzebruch surface with $a\geq 0$, we
let $E$ be a $(-a)$--curve on $\FF_a$ and $F$ a fibre of the
ruling on $\p^1$, so that $F^2=0$ and $E\cdot F=1$. Then a pair
$(X,D)$ with $X=\FF_a$ and $D\equiv \alpha E+\beta F$ will be
called a $(a,\alpha,\beta)$--{\it pair} or an
$(\alpha,\beta)$--{\it pair} on $\FF_a$.\par

Consider a pair $(X,D)$ as above. Let $x_1,...,x_n$ be distinct
points on $X$. Consider the blow--up $p:{\rm
Bl}_{x_1,...,x_n}(X)\to X$ at the given points. On ${\rm
Bl}_{x_1,...,x_n}(X)$ we have the exceptional divisors
$E_1,...,E_n$ corresponding to $x_1,...,x_n$. Consider the divisor
$D_{x_1,...,x_n}:=p^*(D)-E_1-...-E_n$. The pair $({\rm
Bl}_{x_1,...,x_n}(X),D_{x_1,...,x_n})$ will be called the {\it
internal projection} of $(X,D)$ from $x_1,...,x_n$.\par

In the same setting, the pair $({\rm Bl}_{x_1,...,x_n}(X),p^*(D))$
will be called a {\it blow--up} of $(X,D)$.\par

Similarly, consider the divisor
$D_{2x_1,...,2x_n}:=p^*(D)-2E_1-...-2E_n$. The pair $({\rm
Bl}_{x_1,...,x_n}(X),D_{2x_1,...,2x_n})$ will be called the {\it
tangential projection} of $(X,D)$ from $x_1,...,x_n$.\par

\section{Degeneration of projections}\label{degproj}

In this section we generalize some of the ideas presented in \S\S
3-4 of \cite {CMR}, to which we will constantly refer. This will
enable us to prove an extension of Theorem 4.1 of \cite {CMR},
which will be useful later.\par

Let $X\subset \Proj^r$ be an irreducible, non--degenerate
projective variety of dimension $n$. We fix $k\geq 1$, we assume
that $X$ is not $k$--defective and that $s^{(k)}(X)=(k+1)n+k$.\par

Let us fix an integer $s$ such that $r-s^{(k)}(X)\leq s\leq
r-s^{(k-1)}(X)-2$, so that $s^{(k-1)}(X)+1\leq r-s-1\leq
s^{(k)}(X)-1$. Let $L\subset \Proj^r$ be a general projective
subspace of dimension $s$ and let us consider the projection
morphism $\pi_L:S^{k-1}(X)\to \p^ {r-s-1}$ of $X$ from $L$. Notice
that, under our assumptions on $s$, one has:

$$\pi_L(S^k(X))=\p^ {r-s-1}, \quad \pi_L(S^{k-1}(X))\subset\p^ {r-s-1}.$$

Let $p_1,...,p_k\in X$ be general points and let $x\in
<p_1,...,p_k>$ be a general point, so that $x\in S^{k-1}(X)$ is a
general point and $T_{S^{k-1}(X),x}=T_{X,p_1,...,p_k}$. We will
now study how the projection $\pi_L: S^{k-1}(X)\to \Proj^ {r-s-1}$
{\it degenerates} when its centre $L$ tends to a general
$s$-dimensional subspace $L_0$ containing $x$, i.e. such that
$L_0\cap S^{k-1}(X)=L_0 \cap T_{X,p_1,...,p_k}=\{x\}$. To be more
precise we want to describe the {\it limit} of a certain {\it
double point scheme} related to $\pi_L$ in such a degeneration.
\par

Let us describe in detail the set up in which we will work. We let
$T$ be a general $\p^{s^{(k-1)}(X)+s+1}$ which is tangent to
$S^{k-1}(X)$ at $x$, i.e. $T$ is a general $\p^{s^{(k-1)}(X)+s+1}$
containing $T_{X,p_1,...,p_k}$. Then we choose a general line
$\ell$ inside $T$ containing $x$, and we also choose $\Sigma$ a
general $\p^{s-1}$ inside $T$. For every $t\in \ell$, we let $L_t$
be the span of $t$ and $\Sigma$. For $t\in \ell$ a general point,
$L_t$ is a general $\p ^s$ in $\p^r$. For a general $t\in \ell$,
we denote by $\pi_t: S^{k-1}(X)\to \p^ {r-s-1}$ the projection
morphism of $S^{k-1}(X)$ from $L_t$. We want to study the limit of
$\pi_t$ when $t$ tends to $x$. We will suppose from now on that
$k\geq 2$, since the case $k=1$ has been considered in
\cite{CMR}.\par

In order to perform our analysis, consider a neighborhood $U$ of
$x$ in $\ell$ such that $\pi_t$ is a morphism for all $t\in
U\setminus \{x\}$. We will fix a local coordinate on $\ell$ so
that $x$ has the coordinate $0$, thus we may identify $U$ with a
disk around $x=0$ in ${\C}$. Consider the products:

$${\cal X}_1=X\times U,\quad {\cal X}_2=S^{k-1}(X)\times U,\quad \p^
{r-s-1}_U= \p^ {r-s-1}\times U$$

The projections $\pi_t$, for $t\in U $, fit together to give a
morphism $\pi_1: {\cal X}_1\to \p^ {r-s-1}_U$ and a rational map
$\pi_2: {\cal X}_2\map \p^ {r-s-1}_U$, which is defined everywhere
except at the pair $(x,x)=(x,0)$. In order to extend it, we have
to blow up ${\cal X}_2$ at $(x,0)$. Let $p: \tilde{\cal X}_2\to
{\cal X}$ be this blow--up and let $Z\simeq \p^{s^{(k-1)}(X)}$ be
the exceptional divisor. Looking at the obvious morphism $\phi:
\tilde{\cal X}_2\to U$, we see that this is a flat family of
varieties over $U$. The fibre over a point $t\in U\setminus\{0\}$
is isomorphic to $S^{k-1}(X)$, whereas the fibre over $t=0$ is of
the form $\tilde S\cup Z$, where $\tilde S\to S^{k-1}(X)$ is the
blow up of $S^{k-1}(X)$ at $x$, and $\tilde{S}\cap Z=E$ is the
exceptional divisor of this blow up, the intersection being
transverse.\par On $\tilde{\cal X}_2$ the projections $\pi_t$, for
$t\in U$, fit together now to give a {\it morphism} $\tilde \pi:
\tilde{{\cal X}}_2\to \p^ {r-s-1}_U$.
\par

By abusing notation, we will denote by $\pi_0$ the restriction of
$\tilde \pi$ to the central fibre $\tilde S\cup Z$. The
restriction of $\pi_0$ to $\tilde S$ is determined by the
projection of $S^{k-1}(X)$ from the subspace $L_0$: notice in fact
that, since $L_0\cap S^{k-1}(X)=L_0 \cap T_{X,p_1,...,p_k}=\{x\}$,
this projection is not defined on $S^{k-1}(X)$ but it is well
defined on $\tilde S$. \par

As for the action of $\pi_0$ on the exceptional divisor $Z$, this
is explained by the following lemma, whose proof is analogous to
the proof of Lemma 3.1 of \cite {CMR}, and therefore we omit it:

\begin{Lemma} In the above setting, $\pi_0$ maps isomorphically $Z$ to
the $s^{(k-1)}(X)$--dimensional linear space $\Theta$ which is the
projection of $T$ from $L_0$. \label{proiezione} \end{Lemma}

Now we consider ${\cal X}_1\times_U\tilde{{\cal X}}_2$, which has
a natural projection map $\psi: {\cal X}_1\times_U\tilde{{\cal
X}}_2\to U$. One has a commutative diagram:

$$
\begin{array}{ccc}
{\cal X}_1\times_U\tilde{{\cal X}}_2 & \stackrel{\overline
\pi}\rightarrow &
\p^ {r-s-1}_U\\
\psi\downarrow &&\downarrow  \\
U &\stackrel{ { id_U } }\rightarrow & U
\end{array}
$$

\noindent where $\overline \pi=\pi\times\tilde\pi$. For the
general $t\in U$, the fibre of $\psi$ over $t$ is $X\times
S^{k-1}(X)$, and the restriction $\overline \pi_t:X\times
S^{k-1}(X)\to  \p^ {r-s-1}$ of $\overline \pi$ to it is nothing
but $\pi_{t|X}\times \pi_{t|S^{k-1}(X)}$. We denote by
$\Delta_t^{(s,k)}$ the double point scheme of $\overline \pi_t$.
Notice that $\dim(\Delta_t^{(s,k)})\geq s^{(k)}(X)+s-r$ and, by
the generality assumptions, we may assume that equality holds for
all $t\not= 0$. Finally consider the flat limit
$\tilde\Delta_0^{(s,k)}$ of $\Delta_t^{(s,k)}$ inside
$\Delta_0^{(s,k)}$. We will call it the {\it limit double point
scheme} of the map $\overline \pi_t$, $t\not=0$. We want to give
some information about it. Notice the following lemma, whose proof
is similar to the one of Lemma 3.2 of \cite {CMR}, and therefore
we omit it:

\begin{Lemma}  In the above setting, every irreducible component of
$\Delta_0^{(s,k)}$ of dimension $s^{(k)}(X)+s-r$ sits in the limit
double point scheme $\tilde\Delta_0^{(s,k)}$. \label{flatlimit}
\end{Lemma}

Let us now denote by:

\begin{itemize}

\item $X_T$ the scheme cut out by $T$ on $X$. $X_T$ is cut out on
$X$ by $r-s^{(k-1)}(X)-s-1$ general hyperplanes tangent to $X$ at
$p_1,...,p_k$. We call $X_T$ a {\it general
\hbox{$(r-s^{(k-1)}(X)-s-1)$-tangent} section} to $X$ at
$p_1,...,p_k$. Remark that each component of $X_T$ has dimension
at least $n-(r-s^{(k-1)}(X)-s-1)=s^{(k)}(X)+s-r$;

\item $Y_T$ the image of $X_T$ via the restriction of $\pi_0$ to
$X$. By Lemma \ref{proiezione}, $Y_T$ sits in $\Theta=\pi_0(Z)$,
which is naturally isomorphic to $Z$. Hence we may consider $Y_T$
as a subscheme of $Z$;

\item  $Z_T\subset X\times Z$ the set of pairs $(x,y)$ with $x\in
X_T$ and $y=\pi_0(x)\in Y_T$. Notice that $Z_T\simeq X_T$;

\item $\Deltaprimo_0^ {(s,k)}$ the double point scheme of the
restriction of $\pi_0$ to $\tilde S\times X$.
\end{itemize}

With this notation, the following lemma is clear (see Lemma 3.3 of
\cite {CMR}):

\begin{Lemma} In the above setting, $\Delta_0^{(s,k)}$ contains as
irreducible components $\Deltaprimo_0^{(s,k)}$ on $X\times\tilde
S$ and $Z_T$ on $X\times Z$. \label{doublepointscheme}
\end{Lemma}

As an immediate consequence of Lemma \ref{flatlimit} and Lemma
\ref{doublepointscheme}, we have the following proposition (see
Proposition 3.4 of \cite {CMR}):

\begin{Proposition} In the above setting, every irreducible component of
$X_T$, off $T_{X,p_1,\ldots,p_k}$, of dimension $s^{(k)}(X)+s-r$
gives rise to an irreducible component of $Z_T$ which is contained
in the limit double point scheme $\tilde\Delta_0^{(s,k)}$.
\label{limitdoublepointscheme}
\end{Proposition}

\begin{Remark}\label{added}{\rm
We notice that the implicit hypothesis "off $T_{X,x}$" has to be
added also in the statement of Proposition 3.4 of \cite{CMR}.
Actually in the applications in \cite{CMR} this hypothesis is
always fulfilled.} \end{Remark}

So far we have essentially extended word by word the contents of
\S 3 of \cite {CMR}. This is not sufficient for our later
applications. Indeed we need a deeper understanding of the
relation between the double points scheme $\Delta_t^{(s,k)}$ and
$(k+1)$--secant $\p^k$'s to $X$ meeting the centre of projection
$L_t$ and related degenerations when $t$ goes to $0$. We will do
this in the following remark.\par

\begin{Remark}\label{remdoublepointscheme}{\rm (i) It is interesting to give
a different geometric interpretation for the general double point
scheme $\Delta_t^{(s,k)}$, for $t\not= 0$. Notice that, by the
generality assumption, $L_t\cap S^k(X)$ is a variety of dimension
$s^k(X)+s-r$, which we can assume to be irreducible as soon as
$s^k(X)+s-r>0$. Take the general point $w$ of it if
$s^k(X)+s-r>0$, or any point of it if $s^k(X)+s-r=0$. Then this is
a general point of $S^{k}(X)$. This means that $w\in
<q_0,...,q_k>$, with $q_0,...,q_k$ general points on $X$. Now, for
each $i=0,...,k$, there is a point $r_i\in <q_0,...,\hat
q_i,...q_k>$ which is collinear with $w$ and $q_i$. Each pair
$(q_i,r_i)$, $i=0,...,k$, is a general point of a component of
$\Delta_t^{(s,k)}$. Conversely the general point of any component
of $\Delta_t^{(s,k)}$ arises in this way.\medskip

(ii) Now we specialize to the case $t=0$. More precisely, consider
$Z_T\subset X\times Z$ and a general point $(p,q)$ on an
irreducible component of it of dimension $s^{(k)}(X)+s-r$, which
therefore sits in the limit double point scheme
$\tilde\Delta_0^{(s,k)}$. Hence there is a $1$--dimensional family
$\{(p_t,q_t)\}_{t\in U}$ of pairs of points such that
$(p_t,q_t)\in \Delta_t^{(s,k)}$ and $p_0=p, q_0=q$. \par

By part (i) of the present remark, we can look at each pair
$(p_t,q_t)$, $t\not=0$, as belonging to a $(k+1)$--secant $\p^k$
to $X$, denoted by $\Pi_t$, forming a flat family $\{\Pi_t\}_{t\in
U\setminus\{0\}}$ and such that $\Pi_t\cap L_t\not=\emptyset$.
Consider then the flat limit $\Pi_0$, for $t=0$, of the family
$\{\Pi_t\}_{t\in U\setminus\{0\}}$. Since $q\in Z$, clearly
$\Pi_0$ contains $x$. Moreover it also contains $p$. This implies
that $\Pi_0$ is the span of $p$ with one of the $k$--secant
$\p^{k-1}$'s to $X$ containing $x\in S^{k-1}(X)$.} \end{Remark}

As an application of the previous remark, we can prove the
following crucial theorem, which extends Theorem 4.1 of
\cite{CMR}:

\begin{Theorem} \label{boh} Let $X\subset \p^r$ be an irreducible,
non--degenerate, projective variety such that
$s^{(k)}(X)=(k+1)n+k$. Then:

$$d_{X,k}\cdot \deg(X_k)\leq \nu_k(X).$$

\noindent In particular:

\begin{itemize}

\item [(i)] if $r\geq (k+1)(n+1)$ and $X$ is not $k$--weakly
defective, then:

$$\deg(X_k)\leq \nu_k(X);$$

\item [(ii)] if $r=(k+1)n+k$ then:

$$d_{X,k}\leq \mu_k(X).$$\end{itemize}

\end{Theorem}

\begin{proof} We let $s=h^{(k)}(X)=r-s^{(k)}(X)$ and we apply Remark \ref{remdoublepointscheme}
to this situation. Then $X_T$ has $d_{X,k}\cdot \deg(X_k)$
isolated points, which give rise to as many flat limits of
$(k+1)$--secant $\p^k$'s to $X$ meeting a general $\p^s$. By the
definition of $\nu_k(X)$ the first assertion follows. Then (i)
follows from Lemma \ref{birational} and (ii) follows by
(\ref{munu}).\end{proof}
\section {Tangent cones to higher secant varieties}\label{tangentcones}

In this section we describe the tangent cone to the variety
$S^k(X)$, at a general point of $S^l(X)$, where $0\leq l<k$, and
$X\subset\p^r$ is an irreducible, projective variety of dimension
$n$. Our result is the following theorem, which can be seen as a
generalization of Terracini's Lemma:

\begin{Theorem} Let $X\subset \p^r$ be an irreducible, non--degenerate,
projective variety and let $l, m\in\mathbb{N}$ be such that
$l+m=k-1$. If $z\in S^l(X)$ is a general point, then the cone
$J(T_{S^l(X),z},S^m(X))$ is an irreducible component of
$(C_{S^k(X),z})_{\red}$. Furthermore one has:

$$\mult_z(S^k(X))\geq \deg (J(T_{S^l(X),z},S^m(X)))\geq
\deg(S^m(X_{l+1})).$$

 \label{conotang} \end{Theorem}

\begin{proof} We assume that $S^l(X)\neq\p^r$, otherwise the assertion is
trivially true.\par

 The scheme $C_{S^k(X),z}$ is of pure dimension $s^{(k)}(X)$. Let now
$w\in S^m(X)$ be a general point. By Terracini's lemma and by the
generality of $z\in S^l(X)$, we get:

$$\dim(J(T_{S^l(X),z},S^m(X)))=\dim(J(T_{S^l(X),z},T_{S^m(X),w}))=$$
$$=\dim(J(S^l(X),S^m(X)))=\dim(S^k(X))=s^{(k)}(X)$$

\noindent Thus, since $J(T_{S^l(X),z},S^m(X))$ is irreducible and
reduced, it suffices to prove the inclusion
$J(T_{S^l(X),z},S^m(X))\subseteq (C_{S^k(X),z})_{\red}$.

Let again $w\in S^m(X)$ be a general point.  We claim that
$w\not\in T_{S^l(X),z}$. Indeed $S^l(X)\neq\p^r$ and by (\ref
{cone}):

$${\rm Vert}(S^l(X)):=\bigcap_{y\in S^l(X)}T_{S^l(X),y}$$

\noindent  is a proper linear subspace of $\p^r$. If the general
point of $S^m(X)$  would be contained in ${\rm Vert}(S^l(X))$,
then $X\subseteq S^m(X)\subseteq {\rm Vert}(S^l(X))$ and $X$ would
be degenerate, contrary to our assumption.

Since $w\not\in T_{S^l(X),z}$, then $z$ is a smooth point of the
cone $J(w,S^l(X))$. We deduce that:

$$<w,T_{S^l(X),z}>=T_{J(w,S^l(X)),z}=C_{J(w,S^l(X)),z}\subseteq
C_{J(S^m(X),S^l(X)),z}=C_{S^k(X),z}.$$

\noindent By the generality of $w\in S^m(X)$ we finally have
$J(T_{S^l(X),z},S^m(X))\subseteq C_{S^k(X),z}$. This proves the
first part of the theorem.\par

To prove the second part, we remark that:

$$\mult_z(S^k(X))=\deg(C_{S^k(X),z})\geq
\deg(J(T_{S^l(X),z},S^m(X))).$$

\noindent Now, if $p_0,...,p_l\in X$ are general points, then
$J(T_{S^l(X),z},S^m(X))$ is the cone with vertex $T_{S^l(X),z}$
over $\tau_{X,p_0,...,p_l}(S^m(X))$, and, by Lemma \ref {projoin}
we  have that $\tau_{X,p_0,...,p_l}(S^m(X))=S^m(X_{l+1})$. Thus
$\deg(J(T_{S^l(X),z},S^m(X)))\geq\deg(S^m(X_{l+1}))$, proving the
assertion.\end{proof}

\section{A lower bound on the degree of secant varieties}\label{minsecdeg}

As we recalled in \S \ref {notation}, the degree $d$ of an
irreducible non--degenerate variety $X\subset\p^r$ verifies the
lower bound

\begin {equation}\label{mindeg} d\geq  \codim(X) +1.\end{equation}

Varieties whose degree is equal to this lower bound are called
varieties of {\it minimal degree}. As well known, they have nice
geometric properties, e.g. they are rational (see \cite {EH}). In
the present section we will prove a lower bound on the degree of
the $k$--secant variety to a variety $X$. This bound generalizes
(\ref{mindeg}) and we will see that varieties $X$ attaining it
have interesting features which resemble the properties of minimal
degree varieties.\par

Before proving the main result of this section, we need a useful
lemma.

\begin{Lemma}\label{p_1} Let $X\subset\p^r$ be an irreducible,
non--degenerate, projective variety and let $k\geq 0$ be an
integer such that $S^k(X)\neq\p^r$. Let $p\in X$ be an arbitrary
point and, by abusing notation, let $X^1$ be the projection of $X$
from $p$. Then one has:\medskip

\begin{itemize}

\item[(i)] $t_{X,p}(S^k(X))=S^k(X^1)$;\par

\item[(ii)] the general point  in $X$ does not belong to ${\rm
Vert}(S^k(X))$; \par

\item[(iii)]  if $p\in X\setminus(X\cap {\rm Vert}(S^k(X))$, in
particular if $p\in X$ is a general point, then
$t_{{X,p}_{|S^k(X)}}$ is generically finite to its image
$S^k(X^1)$ and $s^{(k)}(X)=s^{(k)}(X^1)$;\par

\item [(iv)] if $X$ is not $k$--defective and $p\in X\setminus
(X\cap {\rm Vert}(S^k(X))$, then $X^1$ is also not
$k$--defective;\par

\item[(v)] if $p\in X\setminus(X\cap {\rm Vert}(S^k(X))$ and if
$\theta_{k}(X)$ denotes the degree of $t_{{X,p}_{|S^k(X)}}$ to its
image, then:

$$\deg(S^k(X))=\theta_{k}(X)\cdot\deg(S^k(X^1))+\mult_{p}(S^k(X))\geq$$
$$\geq \deg(S^k(X^1))+\mult_{p}(S^k(X))$$

\noindent and

$$\mu_k(X^1)=\theta_{k}(X)\cdot \mu_k(X).$$

\end{itemize} \medskip

In particular:\medskip

\begin{itemize}

\item[(vi)] if $p\in X\setminus (X\cap {\rm Vert}(S^k(X))$ and if

$$\deg(S^k(X))=
\deg(S^k(X^1))+\mult_{p}(S^k(X))$$

\noindent then $\theta_{k}(X)=1$, i.e.
$t_{{X,p}_{|S^k(X)}}:S^k(X)\map S^k(X^1)$ is birational and then
$\mu_k(X^1)=\mu_k(X)$;\par

\item[(vii)] if, in addition, $\mu_k(X^1)=1$ then also
$\mu_k(X)=1$ and $\theta_{k}(X)=1$.

\end{itemize} \end{Lemma}

\begin{proof} Part (i) follows by Lemma \ref {projoin}. \par

Since $S^k(X)$ is a proper subvariety in $\p^r$, then ${\rm
Vert}(S^k(X))$ is a proper linear subspace of $\p^r$. This implies
part (ii). Part (iii) is immediate. \par

If $S^k(X)=\p^r$ the assertion is clear. If $X$ is not
$k$--defective and $S^k(X)\not=\p^r$, we have
$s^{(k)}(X)=(k+1)n+k<r$. By part (iii) we have also
$s^{(k)}(X^1)=(k+1)n+k\leq r-1$, i.e. $X^1$ is also not
$k$--defective. This proves (iv).\par

The first assertion of part (v) is immediate.  Furthermore we have
a commutative diagram of rational maps:

$$
\begin{array}{ccc}
S^k_X & \stackrel{t}\dashrightarrow & S^k_{X^1} \\
p^k_X\downarrow &  &\downarrow p^k_{X^1} \\
S^k(X) &\stackrel{ { t_{X,p|S^k(X)} } }\dashrightarrow & S^k(X^1)
\end{array}
$$

\noindent where $t$ is determined, in an obvious way, by
$t_{X,p}$. By the hypothesis, $t_{X,p|S^k(X)}$ has degree
$\theta_{k}(X)$, whereas $t$ is easily seen to be birational.
Hence the conclusion follows. Parts (vi) and (vii) are now
obvious.\end{proof}

Now we come to the main result of this section:

\begin{Theorem} \label{maxbound}
Let $X\subset \p^{r}$, $r=s^{(k)}(X)+h$, $h:=h^{(k)}(X)> 0$, be an
irreducible, non--degenerate, projective variety. Then:

\begin{equation} \label {goodbound}\deg(S^k(X))\geq
\binom{h+k+1}{k+1}\end{equation}

\noindent and, if $l=0,...,k$ and $x\in S^l(X)$ is any point,
then:

\begin{equation} \label {goodbound2} \mult_x(S^k(X))\geq
\binom{h+k-l}{k-l}.\end{equation}

\noindent Suppose equality holds in (\ref {goodbound}) and $h\geq
1$. Then:

\begin{itemize}

\item[(i)] if $x\in X$ is a general point, one has:

$$C_{S^k(X),x}=J(T_x(X),S^{k-1}(X)),\quad
\mult_x(S^{k}(X))=\binom{k+h}{k};$$

\item[(ii)] for every $m$ such that  $1\leq m\leq h$, one has:

$$\deg(S^k(X^m))=\binom{h-m+k+1}{k+1};$$

\item[(iii)]  for every $m$ such that  $1\leq m\leq h$,  the
projection from a general point $x\in X^{m-1}$:

$$t_{{X^{m-1},x}_{|S^k(X^{m-1})}}:S^k(X^{m-1})\map S^k(X^{m})$$

\noindent is birational;

\item[(iv)] for every $m$ such that  $1\leq m\leq k$ one has:

$$\deg(S^{k-m}(X_m))=\binom{h+k-m+1}{k-m+1};$$

\noindent in particular $X_k$ is a variety of minimal degree;

\item[(v)] if $X$ is not $k$--defective, then, for every $m$ such
that  $1\leq m\leq h$, also $X^m$ is not $k$--defective and
$\mu_k(X)=\mu_k(X^m)$;

\item[(vi)] if $X$ is not $k$--defective then:

$$d_{X,k}\leq\mu_k(X).$$
\end{itemize}
\end{Theorem}

\begin{proof}  We make induction on both $k$ and $h$. For $k=0$
we have the bound \ref {mindeg} for the minimal degree of an
algebraic variety, while for $h=0$ the assertion is obvious for
every $k$. Let us project $X$ and $S^k(X)$ from a general point
$x\in X$. By Lemma \ref {p_1}, Theorem \ref{conotang}, Lemma \ref
{kl} and by induction we get:

$$\deg(S^k(X))\geq \deg(S^k(X^1))+\mult_x(S^k(X))\geq$$
$$\geq \deg(S^k(X^1))+\deg(S^{k-1}(X_1))\geq
\binom{k+h}{k+1}+\binom{k+h}{k}=\binom{k+h+1}{k+1}$$

\noindent whence (\ref {goodbound}) follows. Let now $x\in S^l(X)$
be a general point, then by Theorem \ref {conotang}, Lemma \ref
{kl} and by (\ref {goodbound}) one has:

$$\mult_x(S^k(X))\geq \deg (S^{k-l-1}(X_{l+1}))\geq \binom{k+h-l}{k-l}$$

\noindent proving (\ref {goodbound2}) in this case. Of course
(\ref {goodbound2}) also holds if $x\in S^l(X)$ is any point.

If equality holds in (\ref {goodbound}), one immediately obtains
assertions (i)--(iv) for $m=1$. By an easy induction one sees that
(i)--(iv) hold in general.\par

Assertion (v) follows by Lemma \ref{p_1}. As for (vi), consider
the following commutative diagram:

$$
\begin{array}{ccc}
X & \stackrel{\tau_{X,k}}\dashrightarrow & X_k \\
t_{X,h}\downarrow & &\downarrow t_{X_k,h}\\
X^h &\stackrel{\tau_{X^h,k}} \dashrightarrow & \p^n.
\end{array}
$$

\noindent Notice that the vertical maps $t_{X,h}, t_{X_k,h}$ are
birational being projections from $h$ general points on a variety
of codimension bigger than $h$. Thus one has:

$$d_{X,k}=d_{X^h,k}.$$

\noindent On the other hand, by Theorem \ref {boh} and Lemma \ref
{p_1} one has:

$$d_{X^h,k}\leq \mu_k(X^h)=\mu_k(X)$$

\noindent which proves the assertion.\end{proof}

\begin{remark} {\rm It is possible to improve the previous result.
For example, using Lemma \ref{p_1}, one sees that (i) holds not
only if $x\in X$ is general, but also if $x$ is any smooth point
of $X$ not lying on ${\rm Vert}(S^k(X))$. Similar improvements can
be found for (ii)--(v). We leave this to the reader, since we are
not going to use it later.} \end{remark}

\begin{defn}\label{msdeg} {\rm Let $X\subset\p^r$ be an irreducible,
non--degenerate, projective variety of dimension $n$. Let $k$ be a
positive integer.\par

Let $k\geq 2$ be an integer. One says that $X$ is {\it
$k$--regular} if it is smooth and if there is no subspace
$\Pi\subset \p^r$ of dimension $k-1$ such that the scheme cut out
by $\Pi$ on $X$ contains a finite subscheme of length $\ell\geq
k+1$. By definition $1$--{\it regularity} coincides with
smoothness. \par

We say that {\it $X$ has minimal $k$--secant degree}, briefly $X$
is an $\mathcal M^k$--variety, if $r=s^{(k)}(X)+h$,
$h:=h^{(k)}(X)>0$, and $\deg(S^k(X))=\binom{h+k+1}{k+1}$ (compare
with Theorem \ref {maxbound}). \par

We say that $X$ is a variety {\it with the minimal number of
apparent $(k+1)$--secant $\p^{k-1}$'s}, briefly $X$ is an
$\mathcal{MA}^{k+1}_{k-1}$--variety, if $s^{(k)}(X)=(k+1)n+k$,
$r=s^{(k)}(X)+h$, $h:=h^{(k)}(X)> 0$, and if
$\nu_k(X)=\binom{h+k+1}{k+1}$ (compare with Theorem \ref{maxbound}
and \ref{munu}). In other words $X$ is an
$\mathcal{MA}^{k+1}_{k-1}$--variety if and only if it is not
$k$--defective, is an $\mathcal M^k$--variety and $\mu_k(X)=1$.
For example, an $\mathcal M^k$--variety which is not $k$--weakly
defective is an $\mathcal{MA}^{k+1}_{k-1}$--variety (see
Proposition \ref {wd}).\par

We say that $X$ is a variety {\it with one apparent
$(k+1)$--secant $\p^{k-1}$}, briefly $X$ is an
$\mathcal{OA}^{k+1}_{k-1}$--variety, if $r=s^{(k)}(X)=(k+1)n+k$
and $\mu_k(X)=1$.}
\end{defn}

The terminology introduced in the previous definition is motivated
by the fact that, for example,
$\mathcal{OA}^{k+1}_{k-1}$--varieties are an extension of {\it
varieties with one apparent double point} or {OADP}--varieties,
classically studied by Severi \cite {Se} (for a modern reference
see \cite {CMR}).

With this definitions in mind, we have:

\begin{Corollary} Let $k$ be a positive integer. Let $X\subset \p^{r}$,
with $r=s^{(k)}(X)+h$, $h:=h^{(k)}(X)\geq 0$, be an irreducible,
non--degenerate, projective variety of dimension $n$. One has:

\begin{itemize}

\item [(i)] if  $X$ is a $\mathcal{M}^{k}$--variety then for every
$m$ such that $1\leq m\leq h$, the variety $X^m$ is again a
$\mathcal{M}^{k}$--variety;\par

\item [(ii)] if  $X$ is a $\mathcal{MA}^{k+1}_{k-1}$-variety then
for every $m$ such that $1\leq m\leq h-1$, the variety $X^m$ is
again a $\mathcal{MA}^{k+1}_{k-1}$-variety and $X^h$ is a
$\mathcal{OA}^{k+1}_{k-1}$--variety;\par

\item [(iii)] if $X$ is either an
$\mathcal{MA}^{k+1}_{k-1}$-variety or an
$\mathcal{OA}^{k+1}_{k-1}$--variety then $\tau_{X,k}:X\map
X_k\subseteq\p^{n+h}$ is birational and $X_k$ is a variety of
dimension $n$ of minimal degree $h+1$. In particular $X$ is a
rational variety and the general member of the movable part of the
linear system of $k$--tangent hyperplane sections is a rational
variety.\end{itemize}

\label{minimal}
\end{Corollary}

\begin{proof} Parts (i) follows by Theorem \ref {maxbound}, part (ii).
Part (ii) follows by Theorem \ref {maxbound}, parts (ii) and (v).
In part (iii), the birationality of $\tau_{X,k}$ follows by
Theorem \ref {boh}, part (ii). The rest of the assertion follows
by Theorem \ref {maxbound}, part (iv).\end{proof}

\begin{Remark}\label{Bronowski} {\rm In the papers
\cite {Br1} \and \cite {Br}, Bronowski considers the case $k=1$,
$h=0$ and the case $k\geq 2, n=2, h=0$. He claims there, without
giving a proof, that the converse of Corollary \ref {minimal}
holds for $h=0$. We will call this the {\it $k$--th Bronowski's
conjecture}, a generalized version of which, for any $h\geq 0$,
can be stated as follows:} Let $X\subset\p^r$ be an irreducible,
non--degenerate, projective variety of dimension $n$. Set
$h:=h^{(k)}(X)$. If $\tau_{X,k}:X\map X_k\subseteq\p^{n+h}$ is
birational and $X_k$ is a variety of dimension $n$ and of minimal
degree $h+1$, then $X$ is either an
$\mathcal{MA}^{k+1}_{k-1}$--variety or an
$\mathcal{OA}^{k+1}_{k-1}$--variety, according to whether  $h$ is
positive or zero. {\rm We call this the {\it $k$--th generalized
Bronowski's conjecture}}.\par

{\rm Even curve case $n=1$ of this conjecture is still open in
general. The results in \cite {CMR}, \cite {Ru}, \cite {Se}, imply
that the above conjecture is true for $X$ smooth if $k=1, h=0$ and
$1\leq n\leq 3$. The general smooth surface case $n=2$, $k\geq 1$,
$h\geq 0$ follows by the results in \S\S   \ref{oasurf} and
\ref{msurf} (see Corollary \ref{Bron}).  This interesting
conjecture is quite open in general.}\end{Remark}

Bronowski's conjecture would, for example, imply that the converse
of part (ii) of Corollary \ref {minimal} holds. The following
result gives partial evidence for this:

\begin{Proposition} Let $k$ be a positive integer. Let $X\subset
\p^{r+1}$, with $r=(k+1)n+k$, be an irreducible, non--degenerate,
not $k$--defective, projective variety of dimension $n$. If the
general internal projection $X^1$ of $X$ is a
$\mathcal{OA}^{k+1}_{k-1}$--variety, then $X$ is a
$\mathcal{MA}^{k+1}_{k-1}$--variety. \label
{invert}\end{Proposition}

\begin{proof} By part (vii) of Lemma \ref {p_1}, we have that $\mu_k(X)=1$
and $\theta_k(X)=1$. Let $d=\deg(S^k(X))$ and let $p\in X$ be a
general point. Then $t_{ {X,p}_{|S^k(X)} }: S^k(X)\map \p^{r}$ is
a birational map and therefore $\mult_x(S^k(X))=d-1$. Let
$p_0,...,p_{k+1}$ be general points of $X$. Since
$S^{k+1}(X)=\p^{r+1}$, then $S^k(X)$ does not contain
$\Pi:=<p_0,...,p_{k+1}>$. Therefore $S^k(X)$ intersects $\Pi$ in a
hypersurface of degree $d$ with multiplicity $d-1$ at
$p_0,...,p_{k+1}$. This implies that $d\leq k+2$. On the other
hand $d\geq k+2$ by Theorem \ref {maxbound}. This proves the
assertion.\end{proof}

It is interesting to remark that the $\mathcal{M}^k$,
$\mathcal{OA}^{k+1}_{k-1}$ and
$\mathcal{MA}^{k+1}_{k-1}$--properties are essentially preserved
under flat limits:

\begin{Proposition} \label{flatlim} Let $X, X'\subset \p^r$ be reduced,
irreducible, non--degenerate, projective varieties of dimension
$n$, such that $s^{(k)}(X)=s^{(k)}(X')$. Suppose that $X'$ is a
flat limit of $X$ and that $X$ is a $\mathcal{M}^k$--variety
[resp. a $\mathcal{OA}^{k+1}_{k-1}$--variety, a
$\mathcal{MA}^{k+1}_{k-1}$--variety]. then $X'$ is also  a
$\mathcal{M}^k$--variety [resp. a
$\mathcal{OA}^{k+1}_{k-1}$--variety, a
$\mathcal{MA}^{k+1}_{k-1}$--variety] and if
$h^{(k)}(X)=h^{(k)}(X')>0$, then $S^k(X')$ is the flat limit of
$S^k(X)$.
\end{Proposition}

\begin{proof} Suppose $X$ is  a
$\mathcal{M}^k$--variety, so that $h^{(k)}(X)=h^{(k)}(X')>0$. Let
$\Sigma$ be the flat limit of $S^k(X)$ when $X$ tends to $X'$. Of
course $S^k(X)$ is an irreducible component of $\Sigma$, thus by
Theorem \ref {maxbound} we have:

$$\binom{k+h+1}{k+1}\leq \deg(S^k(X'))\leq \deg(\Sigma)=\deg(S^k(X))=
\binom{k+h+1}{k+1}$$

\noindent and therefore the equality has to hold, proving the
assertion.\par

Suppose then $X$ is a $\mathcal{MA}^{k+1}_{k-1}$--variety. The
above argument proves that $S^k(X')$ is the flat limit of
$S^k(X)$. Hence $\mu_k(X')\leq \mu_k(X)=1$, proving that also
$\mu_k(X')=1$, namely the assertion.

The case in which $X$ is a $\mathcal{OA}^{k+1}_{k-1}$--variety is
similar and can be left to the reader.\end{proof}

Finally we point out the following:

\begin {Proposition}\label{linearnorm} Let $X\subset \p^r$ be a variety
with $\mu_k(X)=1$, which is $k$--regular and not $k$--defective.
Then $X$ is linearly normal.\end{Proposition}

\begin{proof} Suppose $X$ is not linearly normal. Then there is a variety
$X'\subset \p^{r+1}$ and a point $p\notin X'$ such that the
projection $\pi$ from $p$ determines an isomorphism  $\pi: X'\to
X$. Now we remark that $p\not\in S^k(X')$ because of the
$k$--regularity assumption on $X$. Furthermore, the assumption
$\mu_k(X)=1$ implies that $\pi: S^k(X')\to S^k(X)$ is also
birational. \par

Set, as usual, $h=h^{(k)}(X)$. Then, by Theorem \ref {maxbound} we
deduce:

$${{k+h+1}\choose {h+1}}=\deg(S^k(X))= \deg(S^k(X'))\geq {{k+h+2}\choose
{h+1}}$$

\noindent a contradiction. \end{proof}

\section {Examples}\label {ex}

In this section we give several examples of
$\mathcal{MA}^{k+1}_{k-1}$ and
$\mathcal{OA}^{k+1}_{k-1}$--varieties.

\begin{example} {\rm {\it Rational normal scrolls}. Let $X=S(a_1,...,a_n)$
be an $n$--dimensional rational normal scroll in $\p^r$. We keep
the notation introduced in subsection \ref {RNS}.\par

We will assume $\sum_{1\leq j\leq n;\; k\leq a_j} (a_j-k)-k-1\geq
0$, otherwise, according to Proposition \ref {defscroll}, one has
$S^k(X)=\p^r$, a case which is trivial for us.

\begin{Claim} If $\sum_{1\leq j\leq n;\;
k\leq a_j} (a_j-k)-k-1\geq 0$, then $X=S(a_1,...,a_n)$ is an
${\mathcal M}^k$--variety.\label{Mkscroll}\end{Claim}

\begin{proof} [Proof of Claim \ref {Mkscroll}] In order to see this, one
may generalize Room's specialization argument (see \cite{Room}, p.
257). Indeed, one has a description of $S^k(X)\subset \p^r$ as a
determinantal variety as follows (see \cite{CJ2}): the homogeneous
ideal of $S^k(X)$ is generated by the minors of order $k+2$ of a
suitable matrix of type $(k+2)\times\sum_{1\leq j\leq n;\; k\leq
a_j} (a_j-k)$ of linear forms, i.e. a suitable {\it Hankel} matrix
of linear forms. Since by Proposition \ref {defscroll} one has
$h:=h^{(k)}(X)=\sum_{1\leq j\leq n;\; k\leq a_j} (a_j-k)-k-1$,
then $S^k(X)$ has, as a determinantal variety, the expected
dimension. Therefore it is a specialization of the variety defined
by the $k+2$ minors of a general matrix of type
$(k+2)\times\sum_{1\leq j\leq n;\; k\leq a_j} (a_j-k)$ of linear
forms, which, as well known (see \cite {acgh}, chapt. II, \S 5),
has degree equal to $\binom{\sum_{1\leq j\leq n;\; k\leq a_j}
(a_j-k)}{k+1}$. As a consequence we have:

$$\deg(S^k(X))=\binom{\sum_{1\leq j\leq n;\; k\leq a_j}
(a_j-k)}{k+1}=\binom{h+k+1}{k+1}$$

\noindent which proves Claim \ref {Mkscroll}. \end{proof}

Next we assume that $X$ is not $k$--defective, i.e., according to
Proposition \ref {defscroll}, that $a_1\geq k$. First we will
consider the case in which $r=(k+1)n+k$, i.e.
$a_1+...+a_n=kn+k+1$, $h:=h^{(k)}(X)=0$, namely $S^k(X)=\p^r$.
Then we make the following:

\begin{Claim}\label{OAscroll} If $a_1\geq k$ and $a_1+...+a_n=kn+k+1$,
then $X=S(a_1,...,a_n)$ is a ${\mathcal
OA}^{k+1}_{k-1}$--variety.\end{Claim}

\begin{proof} [Proof of Claim \ref {OAscroll}] What we have to prove is
that $\mu_k(X)=1$, i.e. that there is a unique $(k+1)$--secant
$\p^k$ to $X$ passing through a general point of $\p^r$.\par

Since $a_1\geq k$, then $|H-kF|$ is generated by global sections
and $h^0(X,{\mathcal
O}_X(H-kF))=\sum_{i=1}^n(a_i+1-k)=k(n+1)+1-n(k-1)=k+n+1$. Let

$$\phi_1=\phi_{|kF|}:X\to\p^{k}=\p(V_1)$$

\noindent and

$$\phi_2=\phi_{|H-kF|}:X\to\p^{k+n}=\p(V_2).$$

\noindent where $V_1=H^0(X,{\mathcal
O}_X(kF))^*,V_2=H^0(X,{\mathcal O}_X(H-kF))^*$. Clearly
$\phi_2(X)=S(a_1-k,\ldots,a_n-k)$, hence $\deg(\phi_2(X))=k+1$.
Let $\phi=\phi_1\times \phi_2$. We get a commutative diagram:

$$
\begin{array}{ccc}
X & \stackrel{\phi}\to & \p^k\times\p^{k+n} \\
\downarrow &  &\downarrow  \\
\p^r& \hookrightarrow & \p^{(k+1)(k+n+1)-1}:=\p_{n,n+k}.
\end{array}
$$

\noindent where the right vertical map is the Segre embedding.\par

Recall that $\p_{n,n+k}=\p(V_1\otimes V_2)=\p({\rm
Hom}(V_1^*,V_2))$. Thus one has a rational map
$\psi:\p_{n,n+k}\map\mathbb{G}(k,n+k)$ which associates to the
class of a rank $k+1$ homomorphism $\xi: V_1^*\to V_2$ the
subspace $\p({\rm Im}(\xi))$ of $\p^{n+k}=\p(V_2)$. \par

One has a natural ${\rm GL}(V_1)={\rm GL}(k+1,\C)$--action on
$V_1\otimes V_2$, which descends to a linear ${\rm
PGL}(k+1,\C)$--action on $\p_{n,n+k}$. From the above description
of the map $\psi$, it is clear that the general fiber of $\psi$ is
a linear space of dimension $k^2+2k$, which is also the closure of
a general orbit of this ${\rm PGL}(k+1,\C)$--action.  More
precisely, if $x\in\p_{k,n+k}$ is a general point, then $x$ is the
class of a homomorphism $\xi: V_1^*\to V_2$, i.e. of a linear
embedding $\iota_{\xi}:\p^k=\p(V_1^*)\to \p^{n+k}=\p(V_2)$. If we
denote by  $\p^k_x$ the image of $\iota_{\xi}$, then the closure
$\Psi_x\simeq \p^{k^2+2k}$ of the fibre of $\psi$ through $x$ can
be interpreted as the linear span of ${\rm
Seg}(k,k)=\p^k\times\p^k_x\subset {\rm Seg}(k,n+k)$. One moment of
reflection shows that this  ${\rm Seg}(k,k)=\p^k\times\p^k_x$ is
an {\it entry locus} in the sense of \cite {Zak}, i.e. it is the
closure of the locus of points of ${\rm Seg}(k,n+k)$ described by
its intersection with the $(k+1)$-secant $\p^k$'s to ${\rm
Seg}(k,n+k)$ passing through $x$.\par

Remark now that $\psi$ is well defined along
$\p^r\subset\p_{k,n+k}$. Indeed, up to projective transformations,
we may assume that $\phi(X)$ contains $k+1$ given general points
of $\p^k\times\p^{k+n}$. Hence we can assume that $\p^r$ contains
an arbitrarily given point of $S^{k}({\rm
Seg}(k,n+k))=\p_{n,n+k}$, e.g. a point where $\psi$ is defined. A
different proof can be obtained as an application of Kempf's
Theorem \ref{Kempf} (see Example \ref{quadricfibrationsextended}
below, we leave the details to the reader). Let us denote by
$\widetilde{\psi}:\p^r\map \mathbb{G}(k,n+k)$ the restriction of
$\psi$ to $\p^r$. \par

We claim that $\widetilde{\psi}$ is dominant. In fact, take $\Pi$
a general $k$--dimensional subspace of $\p^{n+k}=\p(V_2)$. Then
$\Pi$ cuts $\phi_2(X)$ at $k+1$ points $p_0,...,p_k$, which, by
the way, can be interpreted as $k+1$ general points of $X$.
Consider the points $q_i:=\phi_2(p_i)\in\p^k=\p(V_1)$,
$i=0,...,k$. Then one has the embedding $\p^k=\p(V_1^*)\to
\Pi\subset \p^{n+k}=\p(V_2)$, which, for every $i=0,...,k$, maps
the hyperplane $<q_0,...,q_{i-1},q_{i+1},...,q_k>$ to the point
$p_i$. As we saw above, the span $\p^k\times \Pi$ is the fibre of
$\psi$ over the point of $\mathbb{G}(k,n+k)$ corresponding to
$\Pi$.  We thus see that it intersects $X \subset \p_{k,n+k}$ at
the points $p_0,...,p_k$.\par

By the theorem of the dimension of the fibers, the general fiber
of $\widetilde{\psi}$ has dimension $k$. Actually its closure is
the intersection of the linear space $\p^r$ with the general fibre
of $\psi$, which is also a linear space of dimension $k^2+2k$.
Hence we see that this intersection is transversal, i.e. the
closure of the general fiber of $\widetilde{\psi}$ is a $\p^k$. By
the previous analysis we see that it is in fact a $(k+1)$--secant
$\p^k$ to $X$ and that the general such $\p^k$ arises in this way.
\par

In conclusion, since the general $(k+1)$--secant $\p^k$ to $X$ is
the fibre of the rational map $\widetilde{\psi}:\p^r\map
\mathbb{G}(k,n+k)$, we see that there is a unique $(k+1)$--secant
$\p^k$ to $X$ passing through the general point of $\p^r$, i.e.
$\mu_k(X)=1$.\end{proof}

Finally we consider the case $a_1\geq k$ and $r>(k+1)n+k$, i.e.
$a_1+...+a_n>kn+k+1$, $h:=h^{(k)}(X)>0$, thus $S^k(X)\not=\p^r$.
In this case we make the:

\begin{Claim}\label{MAscroll}  If $a_1\geq k$ and $a_1+...+a_n>kn+k+1$,
then $X=S(a_1,...,a_n)$ is a ${\mathcal
MA}^{k+1}_{k-1}$--variety.\end{Claim}

\begin{proof} [Proof of Claim \ref {MAscroll}] Since $X$ is not defective,
by Claim \ref {Mkscroll} all what we have to prove is that
$\mu_k(X)=1$. This easily follows by Lemma \ref {p_1}, (vii), and
Claim \ref {OAscroll}, by making a sequence of general internal
projections. \end{proof} }\label{rationalscrolls} \end{example}

\begin{example}{\rm {\it 2--Veronese fibrations of dimension $n$ and
their internal projections from $h$ points, $1\leq h\leq n+1$}.
Consider $\p(a_1,...,a_n)$, with $0\leq a_1\leq\ldots\leq a_n$ and
$\sum_{i=1}^n a_i\geq 2$. Set $k+1=\sum_{i=1}^na_i+n$ and consider
the map:

$$\phi_1:=\phi_{|H|}:\p(a_1,...,a_n)\to S(a_1,\ldots,a_n)\subset\p^{k}.$$

\noindent Notice that, since $n\leq k-1$, one has
$S(a_1,\ldots,a_n)\neq\p^k$. Furthermore $|H+F|$ is very ample on
$\p(a_1,...,a_n)$ and we can consider the embedding:

$$\phi_2:=\phi_{|H+F|}: \p(a_1,...,a_n)\to
S(a_1+1,\ldots,a_n+1)\subset\p^{k+n}$$

Finally let:

$$\phi_3:=\phi_{|2H+F|}:\p(a_1,...,a_n)\to \p^{r},$$

\noindent where:

$$r=h^0(\p(a_1,...,a_n),({\rm Sym}^2({\cal
O}_{\p^1}(a_1)\oplus...\oplus{\cal O}_{\p^1}(a_n))\otimes {\cal
O}_{\p^1}(1))-1=$$

$$=(n+1)\sum_{i=1}^na_i+n(n+1)=(n+1)(k+1)-1=(k+1)n+k.$$

We set $\phi_3(\p(a_1,...,a_n))= X_{(a_1,...,a_n)}$.

\begin{Claim}\label{duever} $X:=X_{(a_1,...,a_n)}$ is a
$\mathcal{OA}^{k+1}_{k-1}$--variety. \end{Claim}

\begin{proof} [Proof of Claim \ref {duever}]  The verification is
conceptually similar to the case of rational normal scrolls we
worked out in the previous example. Indeed we have a diagram:

$$
\begin{array}{ccc}
\p(a_1,...,a_n) & \stackrel{\phi=\phi_1\times \phi_2}\to &
\p^k\times\p^{k+n} \\
\downarrow &  &\downarrow  \\
\p^{nk+n+k}& \hookrightarrow & \p^{(k+1)(k+n+1)-1}:=\p_{k,k+n}.
\end{array}
$$

Consider the restriction $\widetilde{\psi}$ to $\p^{nk+n+k}$ of
the rational map $\psi: \p_{k,k+n}\map\mathbb{G}(k,k+n)$.

Let us apply Kempf's Theorem \ref {Kempf} to the vector spaces
$V_1=H^0(X,\O_X(H))$, $V_2=H^0(X,\O_X(H+F))$ and
$V_3=H^0(X,\O_X(2H+F))$, where the pairing $V_1\otimes V_2\to V_3$
is the obvious multiplication map. By interpreting the elements of
$V_1, V_2, V_3$ as sections of vector bundles on $\p^1$, one
immediately sees that the pairing is $1$--generic and surjective:
we leave the details to the reader. Then the linear span of
$\phi(X)$ under the Segre embedding is $\p(V_3^*)$. Moreover the
intersection scheme of $S^{k-1}({\rm Seg}(k,{k+n}))=S^{k-1}({\rm
Seg}(\p^k,\p^{k+n}))$ and $\p^{k(n+1)+k}=\p(V_3^*)$ is
irreducible, reduced, of codimension $n+1$ and of degree
$\binom{k+n+1}{k}$  in $\p^{kn+n+k}$. \par

In particular the restriction of $\psi$ is well defined on
$\p^{kn+n+k}$. Then one sees that $S^k(X)=\p^{nk+n+k}$ and
$\mu_k(X)=1$ because the general fiber of $\widetilde{\psi}$ is a
general $(k+1)$--secant $\p^k$ of $X$.\end{proof}

Actually we can prove more:

\begin{Claim} \label{duever2} One has:

\begin{itemize}

\item [(i)] $X:=X_{(a_1,...,a_n)}$ is an
$\mathcal{MA}^k_{k-2}$--variety;

\item [(ii)] the internal projection $X^k$ of $X$ from $h$ points,
$1\leq h\leq n$, is an $\mathcal{MA}^k_{k-2}$--variety;

\item [(iii)] the internal projection of $X^{n+1}$ of $X$ from
$n+1$ points is an $\mathcal{OA}^k_{k-2}$--variety.
\end{itemize} \end{Claim}

\begin{proof} [Proof of Claim \ref {duever2}] By Corollary \ref{minimal},
we need to prove only part (i). For this it suffices to observe
that, as a consequence of the proof of Claim \ref {duever}, one
has that $S^{k-1}(X)$ is a subscheme of the intersection scheme of
$S^{k-1}(\p^k,\p^{k+n})$ and of $\p^{kn+n+k}$. Since these two
schemes are reduced, irreducible and of the same dimension, they
coincide. This yields the desired result:

$$\deg(S^{k-1}(X))=\binom{k+n+1}{k}=\binom{k-1+h^{(k-1)}(X)+1}{k-1+1}.$$

\end{proof}

We notice that, for $n=2$, we have conic bundles. Actually
$\p(a_1,a_2)\simeq \FF_a$, where $a=a_2-a_1$, and $H=E+a_2F$. Then
$2H+F\equiv 2E+(2a_2+1)F=2E+(a+k)F$, where $E$ is a $(-a)$--curve
and $F$ is a ruling, so that $a+k\equiv 1$ (mod.
2).\medskip}\label{quadricfibrationsextended} \end{example}

\begin{example}{\rm {\it 5--Veronese embedding of $\p^2$ and its
tangential projections}. In this example we show that the
5--Veronese embedding $X:=V_{2,5}\subset\p^{20}$ of $\p^2$ and its
general $i$--tangential projections $X_i\subset\p^{20-3i}$, are
smooth $\mathcal{OA}^{k+1}_{k-1}$--surfaces, with $k=6-i$, for
$0\leq i\leq 3$. Notice that $X_3$ is nothing else than the
general $3$--internal projection of $V_{2,4}\subset\p^{14}$, the
$4$--Veronese embedding of $\p^2$.

We will proceed as in the previous examples and we will slightly
modify and adapt to our needs a construction of N. Shepherd-Barron
\cite {SB}. Let us first consider the case of $X=V_{2,5}$. Let us
consider the incidence correspondence

$$F=\{(x,l)\in\p^2\times\p^{2*}\;\; :\;\; x\in l\}.$$

\noindent Then $F$, as a divisor in $\p^2\times \p^{2*}$ sits in
$|\O_{\p^2\times{\p^2}^*}(1,1)|$. Let $p_1$ and $p_2$ denote the
projections of $\p^2\times \p^{2*}$ to the two factors. We will
use the same symbols to denote the restrictions of $p_1$ and $p_2$
to $F$. Let $\phi=\phi_{|\O_F(1,2)|}:F\hookrightarrow\p^{14}$.
Since every fiber of $p_2:F\to\p^{2*}$ is embedded as a line in
$\p^{14}$, we get a morphism $\p^{2*}\to \mathbb{G}(1,14)$, which
is ${\rm PGL}(3,\C)$--equivariant by the obvious action of ${\rm
PGL}(3,\C)$ on $\p^2\times{\p^2}^*$, on $F$, etc. (see \cite{SB}),
and therefore it is an isomorphism to the image. By embedding
$\mathbb{G}(1,14)$ into $\p^{104}$ via the Pl\"ucker embedding,
one has a map $\psi:\p^{2*}\to \p^{104}$, which is an isomorphism
to its image $X$.

\begin{Claim}\label{fiveronese} The image of $\psi$ lands in a $\p^{20}$
and $\psi$ is the $5$--Veronese embedding of $\p^{2*}$ to
$\p^{20}$.\end{Claim}

\begin{proof} [Proof of Claim \ref{fiveronese}]  First of all we notice
that $\psi$ is given by a complete linear system, because it is
clearly ${\rm PGL}(3,\C)$--equivariant. Thus, to prove the claim,
it suffices to show that $\deg(X)=25$. This can be proved by a
direct computation, which we leave to the reader, proving that
$\psi$ is defined by polynomials of degree $5$. However we
indicate here a more conceptual argument (see \cite {SB}, p.
74).\par

Let us introduce the following Schubert cycles in
$\mathbb{G}=\mathbb{G}(1,r)$:

 $$A=\{l\in \mathbb{G}\; :\; l \quad{\rm lies\;\ in\;\ a\;\ given\;\
hyperplane } \},$$

 $$B=\{l\in \mathbb{G}\; :\; l \quad {\rm meets\;\ a\;\ given\;\ linear\;\
space\;\ of \;\ codimension\;\ 3}\},$$

 $$C=\{l\in \mathbb{G}\; :\; l \quad {\rm meets\;\ a\;\ given\;\
linear\;\ space\;\ of \;\ codimension\;\ 2}\}.$$\medskip

\noindent  Then $C$ is a hyperplane section of $\mathbb{G}$ in its
Pl\" ucker embedding and $C^2\sim A+B$. Note that, in our case
$r=14$, we have $\deg(X)=X\cdot C^2=X\cdot A+X\cdot B$.

Notice that:

$$X\cdot
B=\deg(F)=(p_1^*\O_{\p^2}(1)+p_2^*\O_{\p^{2*}}(2))^3=18.$$

Let $H\subset \p^{14}$ be a general hyperplane and let $S=F\cap
H$. Then $S$ is the complete intersection of two divisors of type
$(1,1)$ and $(1,2)$ on $\p^2\times \p^{2*}$. By adjunction $K_S$
is the restriction to $S$ of a divisor of type $(-1,0)$, hence
$K_S^2=2$. Now, $X\cdot A$ is equal to the number of fibers of
$p_2$ lying in $H$, i.e. the number of exceptional curves
contracted by the birational morphism $p_2:S\to \p^{2*}$. Then
$X\cdot A=9-K_S^2=7$.

In conclusion $\deg(X)=18+7=25$ proving Claim \ref{fiveronese}.
\end{proof}

Let us recall now that given a vector space $W$ of odd dimension
$2k+1$, there is a natural rational map $\psi:\p(\Lambda^2W)\map
\p(W^*)$, associating to a general alternating 2--form on $W^*$
its kernel. Then the general fiber of $\psi$ is a linear space and
the map is defined by forms of degree $k$ vanishing to the order
al least $k-1$ along $\mathbb{G}(1,2k)\subset\p(\Lambda^2W)$.\par

Now we are ready to prove the:

\begin {Claim}\label{oaveronese} $X:=V_{2,5}\subset\p^{20}$ is a
$\mathcal{OA}^{7}_{5}$--surface.\end {Claim}

\begin{proof} [Proof of Claim \ref {oaveronese}]  Apply the above remark
to $W=H^0(\O_F(1,2))$, in  order to get a rational map
$\psi:\p^{104}\map\p^{14}$. In \cite{SB}, Lemma 12, it is shown
that the locus of indetermination of $\psi$ does not contain
$S^6(X)=<X>$ (as for the last equality see \cite {CC}, Theorem 1.3
or Example \ref {defectivesurf} below). Thus one has a well
defined rational map $\tilde \psi:<X>=\p^{20}\map\p^{14}$, and
Lemma 13 of \cite {SB} ensures that $\tilde \psi$ is dominant.
Notice that this perfectly fits with the geometry of the
situation. Indeed the closure of a general fiber of $\psi$ is a
$\p^{90}$, cutting $<X>=\p^{20}$ in a linear space of dimension
$90+20-104=6$, which is the general fibre of $\tilde \psi$. On the
other hand, since $\tilde\psi$ is defined by forms of degree $7$
vanishing to the order at least $6$ along $X$, then $\tilde\psi$
contracts every $7$--secant $\p^6$ to $X$. Thus a general
$7$--secant $\p^6$ to $X$ is a general fibre of $\tilde\psi$,
which implies $\nu_6(X)=1$.\end{proof}

We can slightly modify the above construction to show that the
general tangential projection $X_i$ is a
$\mathcal{OA}^{7-i}_{5-i}$--surface, for $i=1,2,3$. We will sketch
the case $i=1$ only, since the others follow by iterating the same
argument.\par

Let $p\in\p^{2*}$ be a general point. We consider the line
$l:=p_2^{-1}(p)$ of $F$. Notice that $p_1(l)$ is the line of
$\p^2$ corresponding to $p$. Consider the projection
$\pi_l:\p^{14}\map\p^{12}$ from $l$ and set $F':=\pi_l(F)$. This
is again a scroll in lines, and the family of lines of $F'$ is
parametrized by a surface $X'\subset\mathbb{G}(1,12)\subset
\p^{77}$.

\begin{Claim}\label{tanprojver} In the above situation, one has that $X'$
is the tangential projection of $X=V_{2,5}$, the $5$--Veronese
embedding of $\p^{2*}$, from the point corresponding to $p$.\end
{Claim}

\begin{proof} [Proof of Claim \ref {tanprojver}]  Set $\p^2\times {\rm
Bl}_p(\p^{2*})\supset\widetilde{F}={\rm Bl}_l(F)\to {\rm
Bl}_p(\p^{2*})$ and let $\phi:\widetilde{F}\to\p^{12}$ be the map
given by the linear system
$|p_1^*(\O_{\p^2}(1))+p_2^*(\O_{\p^{2*}}(2))-\widetilde{E}|$,
where $p_1$ and $p_2$ are the projections of $\p^2\times {\rm
Bl}_p(\p^{2*})$ and $\widetilde{E}$ is the exceptional divisor of
$\widetilde{F}$. Then $F'\simeq \phi(\widetilde{F})$ from which it
follows that $X'\simeq {\rm Bl}_p(\p^{2*})$.

Now, the map $\pi_l:\p^{14}\map\p^{12}$ gives rise to a map
$\widetilde{\pi_l}:\mathbb{G}(1,14)\map\mathbb{G}(1,12)$ which is
nothing but the tangential projection of $\mathbb{G}(1,14)$ from
the point corresponding to $l$. This implies that the inclusion
$X'\subset\mathbb{G}(1,12)\subset \p^{77}$ is given by the
pull--back on $X'$ of a linear system of quintics of $\p^{2*}$
which are singular at $p$. To prove the claim it suffices to
remark that the embedding $X'\subset\mathbb{G}(1,12)\subset
\p^{77}$ is given, as usual, by a complete linear system. Moreover
one has $\deg(X')=21$. To see this we have to make exactly the
same calculation as for the computation of $\deg(X)$. In the
present case one has that $X'\cdot B=\deg(F')=15$ and $X'\cdot
A=6$ so that $\deg(X')=21$.

Now we notice that $<X_1>=\p^{17}=S^5(X_1)$ (use Terracini's lemma
or Theorem 1.3 from \cite {CC} or Example \ref {defectivesurf}
below). Arguing as for $X$, we have now a map
$\psi:\p^{77}\map\p^{12}$ which is defined by forms of degree 6
vanishing to the order 5 along $\mathbb{G}(1,12)$. One proves that
$<X_1>$ does not lie in the indeterminacy locus of $\psi$ so that
one has a well defined rational map $\tilde\psi:
<X_1>=\p^{17}\map\p^{12}$ and one shows that this map is dominant.
The fibres of $\tilde\psi$ are the $6$--secant $\p^5$'s to $X_1$,
and therefore $\nu_5(X_1)=1$. \end{proof}}\label{quintics}
\end{example}

\begin{example}{\rm {\it 4--Veronese embedding of $\p^2$ and its
internal projections}. In this example we note that $V_{2,4}$ is a
${\mathcal {MA}}^{4}_2$--surface. This can be proved by using the
formulas in \cite {ES} and \cite {Le} to prove that
$\deg(S^3(V_{2,4}))=35$. By Theorem \ref {maxbound}, (ii), we see
that also that a general $i$--internal projection of $V_{2,4}$,
$i=1,2$, has the same property.\par

Another interesting property of $V_{2,4}$ is that it is
$4$--defective and $S^4(V_{2,4})$ is a hypersurface in $\p^{14}$
(see \cite {CC}, Theorem 1.3 or Example \ref {defectivesurf}
below). One has $\deg(S^4(V_{2,4}))=6$, hence $V_{2,4}$ is a
$\mathcal M^4$--surface. This can be proved as follows. Look at
$V_{2,4}$ as that $2$--Veronese embedding of $V_{2,2}\subset
\p^5$. Thus $S^4(V_{2,4})\subseteq <V_{2,4}>\cap S^4(V_{5,2})$,
where $S^4(V_{5,2})$ is a hypersurface of degree $6$. Notice that
$<V_{2,4}>$ is not contained in $S^4(V_{5,2})$. In fact, since
$V_{2,2}$ is non--degenerate in $\p^5$, then given $6$ general
points of $V_{5,2}$ we can suppose that $V_{2,4}$ contains them.
Thus  we may assume that $<V_{2,4}>$ contains a general point of
$S^5(V_{5,2})=\p^{20}$ which can be chosen to be off
$S^4(V_{5,2})$. Finally we know,
 by Theorem \ref {maxbound}, that $\deg(S^4(V_{2,4}))\geq 6$. This implies
that $S^4(V_{2,4})$ is the scheme--theoretic intersection of
$<V_{2,4}>$ and $S^4(V_{5,2})$ and that $\deg(S^4(V_{2,4}))=6$.
\par

Using this same line of argument, one can give a direct, more
geometric proof that $\deg(S^3(V_{2,4}))=35$. We leave the details
to the reader. }\label{quartics}\end{example}

\begin{example} {\rm {\it The 3--Veronese embedding of the quadric
surface in $\p^3$.} Let $X\subset\p^{15}$ be the $3$--Veronese
embedding of a smooth quadric surface $Q\subset\p^3$. Then $X$ is
a $\mathcal{MA}^5_3$--surface, i.e. $S^4(X)\subset\p^{15}$ is a
hypersurface of degree 6. Indeed, the projection of $X$ from a
point on it is isomorphic to the 2--tangential projection of the
$5$-Veronese embedding of $\p^2$, which is a $\mathcal
{OA}^5_3$--surface, see Example \ref{quintics}. The conclusion
follows from Proposition \ref{invert}.\par

By applying Proposition \ref {flatlim}, one sees that also the
$3$--Veronese embedding of a quadric cone in $\p^3$ is a
$\mathcal{MA}^5_3$--surface.}\label{treverquadrics}
\end{example}

\begin{example}{\it Defective surfaces}. {\rm  The fact that $V_{2,4}$ is
a $\mathcal M^4$--surface is a particular case of a more general
family of examples of surfaces with minimal secant degree.\par

According to Theorem 1.3 from \cite {CC}, this is the list of
$k$--defective surfaces $X\subset \p^r$:

\begin{itemize}

\item [(i)] $r=3k+2$ and $X$ is the $2$--Veronese embedding of a
surface of degree $k$ in $\p^{k+1}$, and $\delta_k(X)=1$;\par

\item [(ii)] $X$ sits in a $(k+1)$--dimensional cone over a curve.
\end{itemize}

We claim that the surfaces of type (i) are $\mathcal
M^k$--surfaces. In fact such a $X$ is contained in $V_{k+1,2}$ and
therefore $S^k(X)\subseteq <X>\cap S^k( V_{k+1,2})$. Here again we
have that:

\begin{itemize}

\item  $<X>$ is not contained in $S^k(V_{k+1,2})$;

\item $S^k(V_{k+1,2})$ is a hypersurface of degree $k+2$, i.e. it
is the set of singular quadrics in $\p^{k+1}$;

\item  $\deg(S^k(X))\geq k+2$, by Theorem \ref {maxbound}.

\end{itemize}

These three facts together imply that the hypersurface $S^k(X)$ is
the scheme--theoretic intersection of $<X>$ and $S^k(V_{k+1,2})$
and that $\deg(S^k(X))=k+2$.
\par

The first instance of this family of examples, obtained for $k=1$,
is the Veronese surface $V_{2,2}$ in $\p^5$, whose secant variety
is a hypersurface of degree 3.}\label{defectivesurf}\end{example}

\begin{example}{\it Weakly defective surfaces}. {\rm  The previous example
can be further extended. \par

According to Theorem 1.3 from \cite {CC}, this is the list of
$k$--weakly defective, not $k$--defective, surfaces $X\subset
\p^r$:

\begin{itemize}

\item [(i)] $r=9$, $k=2$ and $X$ is the $2$--Veronese embedding of
a surface of degree $d\geq 3$ in $\p^{3}$;\par

\item [(ii)]  $r=3k+3$, and $X$ is the cone over a $k$--defective
surface of type (i) in Example \ref {defectivesurf};\par

\item [(iii)] $r=3k+3$, and $X$ is the $2$--Veronese embedding of
a surface of degree $k+1$ in $\p^{k+1}$;\par

\item [(iv)] $X$ sits in a $(k+2)$--dimensional cone over a curve
$C$, with a vertex of dimension $k$. \end {itemize}

We claim that the surfaces of types (i), (ii) and (iii) are
$\mathcal M^k$--surfaces.\par

If $X$ is a surface of type (i), one immediately sees that
$S^2(X)=S^2(V_{3,2})$, hence $\deg(S^2(X))=4$ and $X$ is therefore
a ${\mathcal M}^2$--surface. \par

If $X$ is a surface of type (ii), then $S^k(X)$ is the cone over
the $k$--secant variety of a $k$--defective surface of type (i) in
Example \ref {defectivesurf}. Hence we have $\deg(S^k(X))=k+2$ and
$X$ is a ${\mathcal M}^k$--surface.\par

If $X$ is a surface of type (iii), the same argument we made in
Example \ref {defectivesurf} proves our claim. We leave the
details to the reader.}\label{wdefectivesurf}\end{example}

\begin{example} {\rm {\it Del Pezzo surfaces}. In this example we remark
that smooth del Pezzo surfaces of degree $r$ in $\p^r$,
$r=5,...,9$, are ${\mathcal MA}^{2}_0$--surfaces. This can be
easily seen by applying the double point formula. Proposition
\ref{flatlim} implies that also singular del Pezzo surfaces are
${\mathcal MA}^{2}_0$--surfaces.\par

The Veronese surface $X:=V_{2,3}$ is also an ${\mathcal
MA}^{3}_1$--surface, as can be seen by applying Le Barz's formula
\cite {LB}. However this is a classical result. Indeed
$S^2(V_{2,3})$ is the hypersurface of $\p^9$ consisting of all
cubics which are sums of three cubes of linear forms. These are
the so--called {\it equihanarmonic} cubics, i.e. those
characterized by the vanishing of the $J$--invariant. It is
classically well known that there are four equihanarmonic cubics
in a general pencil (see \cite {EC}, p. 194), i.e. $\deg
(S^2(V_{2,3}))=4$, which means that $V_{2,3}$ is a ${\mathcal
MA}^{3}_1$--surface.\par

We can also give a more geometric proof of this fact by applying
the ideas we have developed so far. Indeed, the general internal
projection $X^1$ of $X$ is the embedding of $\FF_1$ in $\p^8$ via
the linear system $|2E+3F|$. This, according to Example
\ref{quadricfibrationsextended},  is a $\mathcal
{OA}^3_1$--surface. Thus $X$ is  a ${\mathcal MA}^{3}_1$--surface
by Proposition \ref {invert}.}\label{delpezzo}\end{example}

\begin{example} Cones. {\rm Let $X\subset\p^r\subset\p^{r+l+1}$, $l\geq
0$, be an irreducible variety of dimension $n$ which is
non--degenerate in $\p^r$. Let $L=\p^l\subset\p^{r+l+1}$ be such
that $L\cap\p^r=\emptyset$. Let $Y=J(L,X)$ be the cone over $X$
with vertex $L$. Then $\dim(Y)=n+l+1$. More generally for every
$k\geq 1$ we have $S^k(Y)=S(L,S^k(X))$ so that
$s^{(k)}(Y)=s^{(k)}(X)+l+1$. Therefore
$h^{(k)}(Y)=r+l+1-s^{(k)}(Y)=r-s^{(k)}(X)=h^{(k)}(X)$. Moreover
$\deg(S^k(Y))=\deg(S^k(X))$ for every $k\geq 1$. In particular $X$
has minimal $k$--secant degree if and only if $Y$ has also minimal
$k$--secant degree.

For instance, a rational normal scroll $X=S(a_1,...,a_n)$ is a
variety of minimal $k$--secant degree if the least positive
integer $a_i$ is greater or equal than $k$ (see Example
\ref{rationalscrolls}). }\label{cones}
\end{example}

The next example is a slight modification of the previous one. It
shows that some of the hypotheses we will make in our
classification theorems in \S\S  \ref {oasurf}, \ref {msurf} are
well motivated. The first instance of this example, i.e. the case
$k=1$, is due to A. Verra, who kindly communicated it to us. It
could be easily generalized to higher dimensions and codimensions:
we leave the details to the reader.

\begin{example}{\rm Let $C\subset\p^{2k+1+h}\subset\p^{3k+2+h}$, $k\geq 1,
h\geq 0$, be an irreducible curve, non--degenerate in
$\p^{2k+1+h}$. Take $\Pi=\p^k\subset\p^{3k+2+h}$ such that
$\Pi\cap\p^{2k+1+h}=\emptyset$ and a  morphism $\phi:C\to
C'\subset\p^k$ and take $X=\cup_{p\in
C}<p,\phi(p)>\subset\p^{3k+2+h}$. Then $\nu_k(X)=\nu_k(C)$. This
is an exercise in projective geometry which we leave to the
reader. In particular, from Example \ref {rationalscrolls} and
from Theorem \ref {curvess} below, we deduce that $\nu_k(X)=1$ if
and only if $C$ is a rational normal curve. As soon as $k\geq 3$,
one can take as $\phi$ a general projection of $C$ and obtain
examples of smooth surfaces $X\subset\p^{3k+2+h}$, which are not
linearly normal. Let us remark that such a surface $X$ is
$k$--weakly defective, being contained in a cone of vertex a
$\p^k$ over the curve $C$, see \cite{CC}, Theorem 1.3 and Example
\ref {wdefectivesurf}.}\label{wdd}
\end{example}

\section{Classification of curves with minimal secant
degree}\label {curves}

In this section we take care of the classification of curves with
minimal $k$--secant degree.\par

Let $C\subset\p^r$ be an irreducible non--degenerate curve. Then
$C$ is never weakly defective (see \cite {CC}), hence it is never
defective, thus $s^{(k)}(C)=\min\{2k+1, r\}$. The classification
of curves with minimal $k$--secant degree is given by the
following:

\begin{Theorem} \label{curvess} Let $C\subset\p^r$ be an irreducible
non--degenerate curve. Let $k\geq 1$ be an integer such that
$2k+1\leq r$. Then $C$ is an $\mathcal{MA}^{k+1}_{k-1}$ or an
$\mathcal{OA}^{k+1}_{k-1}$--variety if and only if $C$ is a
rational normal curve.
 \end{Theorem}
\begin{proof} As we saw in Example \ref{rationalscrolls}, a
rational normal curve is an $\mathcal{MA}^{k+1}_{k-1}$ or an
$\mathcal{OA}^{k+1}_{k-1}$--variety. \par

Suppose, conversely, that $C$ is an $\mathcal{MA}^{k+1}_{k-1}$ or
an $\mathcal{OA}^{k+1}_{k-1}$--variety. In the latter case, i.e.
if $r=2k+1$, then the assertion is Theorem 3.4 of \cite {CJ}. In
the former case, i.e. if  $h=r-2k-1>0$, then part (ii) of
Corollary \ref{minimal} tells us that $C^h$ is  an
$\mathcal{OA}^{k+1}_{k-1}$--variety. Since, as we saw, $C^h$ is a
rational normal curve, then $C$ itself is a rational normal curve,
proving the assertion.\end{proof}

\begin {Remark} \label{curveregularity} {\rm Notice that, in the
hypotheses of Theorem \ref{curvess}, the rationality of $C$
follows by Corollary \ref{minimal}. If one adds the hypothesis
that $C$ is $k$--regular, then the assertion follows right away
from Proposition \ref{linearnorm}.}\end{Remark}

\section{On a theorem of Castelnuovo--Enriques} \label{CastEnr}

The next sections will be devoted to the classification of
$\mathcal{OA}^{k+1}_{k-1}$--surfaces and
$\mathcal{M}^{k}$--surfaces. For this we will need some
preliminaries, which we believe to be of independent interest,
concerning linear systems of curves on a surface. Indeed the
present section is devoted to review, and improve on, a classical
theorem of Enriques, which in turn generalizes to arbitrary
surfaces an earlier result proved by Castelnuovo for rational
surfaces, see \cite{En} and \cite{Ca}. The expert reader will find
relations between the results of this section and the ones in
\cite{hart} and in \cite{Re}. We will freely use here the notation
introduced in subsections \ref{gen} and \ref{gen'}.\par

The basic tool in this section is Proposition  \ref{adj} below.
This result essentially goes back to Iitaka \cite{Iitaka} and
Dicks \cite{Dicks}, Theorem 3.1, though under the stronger
assumption that  $D$ is an irreducible smooth curve. The case $D$
ample is also well known in the literature, e.g. see  \cite{io2}.
The short proof below, based on Mori's theory, is essentially the
same as in \cite{Dicks}, and we included it here for the reader's
convenience.

\begin{Proposition} \label{adj} Let $X$ be a smooth, irreducible,
projective surface. Let $D$ be a nef divisor on $X$. Set $d:=D^2$,
$g:=p_a(D)$. Assume the pair $(X,D)$ is minimal, not a $h$--scroll
with $h\leq 1$ and it is not a $m$--Veronese pair with $m\leq 2$.
Then $K+D$ is nef and therefore:

\begin{itemize}

\item[(i)] $d\leq 4(g-1)+K^2$;

\item[(ii)] $g\geq 1$ and equality holds if and only if $K$ and
$D$ are numerically dependent and either $d=0$ or $(X,D)$ is a del
Pezzo pair.

\end{itemize}
\end{Proposition}

\begin{proof} Let $C$ be a curve on $X$ such that $C\cdot (K+D)<0$. Since
$D$ is nef, one has $K\cdot C<0$. By Mori's cone Theorem (see
\cite {Mori}, Theorem 1.4), the curve $C$  is a linear combination
of extremal rays. More precisely, there are extremal rays
$E_1,...,E_h$ such that $C\sim\sum_{i=1}^hm_iE_i$, with
$m_1,...,m_h$ positive real numbers. Thus there is one of the
extremal rays $E_1,...,E_h$, e.g. $E:=E_1$ such that $E\cdot
(K+D)<0$. Now one concludes by separately discussing the various
possibilities for $E$ (cfr.  \cite {Mori}, Theorem 2.1):

\begin{itemize}

\item if $E$ is a $(-1)$--curve, one has $K\cdot E=-1$ and
therefore $D\cdot E=0$, against the minimality of $(X,D)$;

\item if $E\simeq \p^1$ and $E^2=0$, one has $K\cdot E=-2$ and
therefore $D\cdot E\leq 1$, against the fact that $(X,D)$ is not a
$h$--scroll for $h\leq 1$;

\item if $E\simeq \p^1$ e $E^2=1$, one has $K\cdot E=-3$ and
therefore $1\leq D\cdot E\leq 2$, against the fact that $(X,D)$ is
not a
 $m$--Veronese with $m\leq 2$.

\end{itemize}

Now notice that:

\begin{equation} \label{quadrato} (K+D)^2=K^2+4(g-1)-d\end{equation}

Since $K+D$ is nef, one has $(K+D)^2\geq 0$, so that
\begin{equation} \label{bound} d\leq 4(g-1)+K^2,\end{equation}
proving (i).

 Similarly,
since $K+D$ is nef, one has $2g-2=(K+D)\cdot D\geq 0$, proving the
first assertion of (ii). If $g=1$, one has $(K+D)\cdot D=0$. Then
the Hodge index theorem implies that $K+D$ and $D$ are numerically
dependent, thus $K\sim lD$, for some rational number $l$. If $d>0$
then $0=(K+D)\cdot D=(l+1)d$ implies $l=-1$ and $(X,D)$ is a del
Pezzo pair. Conversely if $(X,D)$ is a del Pezzo pair then $g=1$.
Similarly, if $d=0$ and $K$ and $D$ are numerically dependent, one
has $g=1$.\end{proof}

\begin{Corollary} \label{gi}  Let $X$ be a smooth, irreducible,
projective surface. Let $D$ be a nef divisor on $X$. Assume the
pair $(X,D)$ is not a $h$--scroll with $h\leq 1$. Set $g:=p_a(D)$.
Then $g\geq 0$ and $g=0$ if and only if $(X,D)$ is obtained by a
$m$--Veronese with $m\leq 2$ with a sequence of blowing--ups.
\end{Corollary}

\begin{proof} By iterated contractions of $(-1)$--curves $E$ such that
$E\cdot D=0$, we arrive to a minimal pair $(X',D')$ such that
$(X,D)$ is obtained from $(X',D')$ with a sequence of
blowing--ups. Moreover $g':=p_a(D')=g$. Notice that $(X',D')$, as
well as $(X,D)$, is not a $h$--scroll with $h\leq 1$. Then the
assertion follows by the second part of Proposition
\ref{adj}\end{proof}

As a consequence we have the following result, essentially due to
Castelnuovo \cite{Ca} and Enriques \cite{En}. The bound
\eqref{grado} was also obtained by Hartshorne, \cite{hart}
Corollary 2.4 and Theorem 3.5, under the assumption $D$ smooth
irreducible curve. Hartshorne does not consider the classification
of the extremal cases, as in \cite{Ca}, but he remarks that the
bound is sharp looking at the cases (i) and (iv) with $a=0$,
Example in {\it loc. cit.}, p. 121. All the results of \cite{hart}
are now straightforward consequences of  Proposition \ref{adj}.

\begin{Theorem} \label{Enr} Let $X$ be a smooth, irreducible,
projective surface. Let $D$ be an irreducible curve on $X$. Set
$d:=D^2$, $g:=p_a(D)$, $r:=\dim(|D|)$. Assume $d\geq 0$ and the
pair $(X,D)$ is not a $h$--scroll with $h\leq 1$. Then:

\begin{equation} \label{grado} d\leq 4g+4+\epsilon
\end{equation}

\noindent where $\epsilon=1$ if $g=1$ and $\epsilon=0$ if $g\neq
1$. Consequently one has:

\begin{equation} \label{erre} r\leq 3g+5+\epsilon
\end{equation}

\noindent and the equality holds in (\ref {grado}) if and only if
it holds in (\ref{erre}).\par

If, in addition, the pair $(X,D)$ is minimal, then the equality
holds in (\ref{grado}), or equivalently in (\ref {erre}), if and
only if one of the following happens:\par

\begin{itemize}

\item [(i)] $g=0$, $r=5$, and $(X,D)$ is a $2$--Veronese pair;\par

\item [(ii)] $g=1$, $r=9$, and $(X,D)$ is a $3$--Veronese
pair;\par

\item [(iii)] $g=3$, $r=14$, and $(X,D)$ is a $4$--Veronese
pair;\par

\item [(iv)] $(X,D)$ is a $(2,a+g+1)$--pair on $X\simeq  \FF_a$,
$a\geq 0$.\end{itemize} \end{Theorem}

\begin{proof}  By arguing as in the proof of
Corollary \ref{gi} we may, and will, assume that the pair $(X,D)$
is minimal. Then note that if $(X,D)$ is a $m$--Veronese with
$m\leq 2$, both \eqref{grado} and \eqref{erre} hold. So we may
assume $(X,D)$ is not a $m$--Veronese with $m\leq 2$.\par

Let us now prove \eqref{grado}. The divisor $D$ is nef so that
bound  \eqref{bound} holds.\par

Assume that $d> 4g+4+\epsilon$. Then $K\cdot D=2g-2-d\leq
-2g-6-\epsilon<0$. Therefore $\kappa(X)=-\infty$. Moreover:

$$4g+4+\epsilon< d\leq 4g-4+K^2$$

yields $K^2\geq 9+\epsilon$. Therefore $\epsilon=0$, i.e. $g\neq
1$, $K_X^2=9$ and $X\simeq \Proj^2$. Hence $D\in |\O_{\Proj
^2}(m)|$, with $m\geq 4$, since $(X,D)$ is not a Veronese pair
with $m\leq 2$ and $g\neq 1$. For such a $D$ one has $m^2=d\leq
4g+4=2m^2-6m+8$. This contradiction proves  (\ref{grado}).\par

Next we remark that \eqref{grado} implies \eqref{erre}. Indeed,
since the general curve $D\in |D|$ is irreducible, by
Riemann--Roch theorem we have $r\leq \max\{d-g+1,g\}$, which
implies (\ref{erre}). \par

Let us prove now that equality holds in \eqref {grado} if and only
if equality holds in \eqref{erre}. The above argument shows that
if equality holds in \eqref{grado} then it holds in \eqref{erre}.
Conversely, if equality holds in (\ref{erre}) then Riemann-Roch
theorem implies that $d-g+1\geq r$ and equality holds in
(\ref{grado}).\par

Finally, suppose equality holds in (\ref{grado}). Then reasoning
as above we deduce $\kappa(X)=-\infty$ and $K^2\geq 8+\epsilon$.
Therefore if  $g=1$ one has $K^2=9$, $(X,D)$ is a del Pezzo pair
and we are in case (ii). We can thus suppose $\epsilon=0$ in
\eqref{grado} and hence $K^2\geq 8$.\par

If $K^2=9$, then $X\simeq \Proj^2$, $D\in |\O_{\Proj ^2}(m)|$,
with $m\geq 1$. The equality $d=4g+4$ is translated into
$m^2=2m^2-6m+8$, so that $m=2$ or $4$ and we get  cases (i) and
(iii).\par

Assume that $K^2=8$. Thus $X\simeq \FF_a$, $a\geq 0$. Furthermore
\eqref{quadrato} shows that  $(K+D)^2=0$ holds. One has:

 $$D\sim \alpha
E+\beta F,$$

\noindent where $E$ is a $(-a)$-curve and $F$ a fibre of the
ruling of $\FF_a$, with $\beta\geq a\alpha$ because $D\cdot E\geq
0$, and $\alpha\geq 2$ since the pair $(X,D)$ is not a scroll. On
the other hand:

$$K\sim-2E-(a+2)F$$

\noindent and therefore:

$$K+D\sim(\alpha-2)E+(\beta-a-2)F.$$

If $\alpha=2$ then adjunction formula implies:

$$\beta=a+g+1$$

\noindent i.e. the assertion. Now:

$$(K+D)^2=(\alpha-2)(2\beta-a\alpha-4).$$

If $a=0$, $(K+D)^2=0$ implies either $\alpha=2$ or $\beta =2$, and
we are done. If $a=1$, the minimality condition yields $\beta\geq
\alpha+1$. Therefore $(K+D)^2=0$ implies $\alpha=2$, and we are
done again. If $a\geq 2$, one has $2\beta-a\alpha-4\geq
a\alpha-4=2(\alpha-2)$. Then $(K+D)^2=0$ implies $\alpha=2$, and
we conclude as above. \end{proof}

\begin{Remark} \label{improve} {\rm Proposition \ref {adj} can be improved.
Indeed, we can prove that if one adds the hypothesis that $D$ is
effective and big, then  $K+D$ is also effective. This can be seen
as a wide extension of the results in \cite {Ba}, p. 196-200.
Following the ideas in \cite {Ca2} one can even give suitable,
interesting lower bounds for $(K+D)^2$. \par

It is also possible to partly extend Proposition \ref{adj} to
higher dimensional varieties. \par

The hypothesis $D$ effective and irreducible in Theorem \ref{Enr}
is essentially used to prove that \eqref{grado} implies
\eqref{erre} and it is too strong. Indeed, we can prove that it
suffices to assume that either $g\not=1$ or $d>0$. However the
proof, based on the aforementioned extensions of Proposition
\ref{adj} as indicated in \cite{Ca2}, is rather long and we
decided not to put it here. We plan to come back to this and to
other extensions of Proposition \ref{adj} and Theorem \ref{Enr} in
the future.}\end{Remark}

\begin{defn} {\rm If the pair $(X,D)$ is as in (iv) of Theorem \ref{Enr},
we will say that it is a $(a,g)$--{\it Castelnuovo pair}
and the corresponding surface $\phi_{|D|}(X)\subset \p^{3g+5}$ of
degree $d=4g+4$, with hyperelliptic hyperplane sections, will be
called an $(a,g)$--{\it Castelnuovo surface} and denoted by
$X_{a,g}$. The motivation for this definition resides in the fact
that Castelnuovo first considered these pairs in his paper \cite
{Ca}. In general, a pair like in (i), (ii), (iii) or (iv) of
Theorem \ref{Enr}, will be called a {\it Castelnuovo extremal
pair}.

We notice that pairs $(X,D)$ as in (ii), (iii) or (iv) can be
characterized as those with $D$ effective, irreducible and nef for
which the hypotheses of  Proposition \ref{adj} are met, so that
$K+D$ is nef, but $K+D$ is not big.}\end{defn}

\begin{Remark} \label{castsurf} {\rm An $(a,k)$--{Castelnuovo surface}
$X_{a,k}$ is $(k+1)$--defective as soon as $a+1+k\equiv 0$ (mod.
2) (see case (i) of Theorem 1.3 of \cite {CC} and Example
\ref{defectivesurf}). In this case the Castelnuovo surface will be
said to be {\it even}. Instead $X_{a,k}$ is an ${\mathcal
OA}^{k+2}_k$ surface if $a+1+k\equiv 1$ (mod. 2), and then the
Castelnuovo surface will be said to be {\it odd}. In fact in this
case $X_{a,k}$ is one of the surfaces described in Example
\ref{quadricfibrationsextended}.

Note that an $(a,k)$--{Castelnuovo surface} $X_{a,k}$ is smooth
unless $k=a-1$, in which case the Castelnuovo surface is even and
it is the $2$--Veronese embedding of a cone over a rational normal
curve of degree $a$.}\end{Remark}

It is useful to point out the following immediate corollaries,
whose easy proofs can be left to the reader:

\begin{Corollary} Let $X$ be a smooth, irreducible,
projective surface. Let $\L$ be a linear system of dimension $r>0$
whose general divisor $D$ is irreducible with geometric genus $g$.
Let $(X',\mathcal L')$ be the resolution on $(X,\mathcal L)$.
Suppose $(X',\mathcal L')$ is not a scroll. Then (\ref {erre}) of
Theorem \ref {Enr} holds. If, in addition, $(X',\mathcal L')$ is
minimal and equality holds in (\ref {erre}), then $(X,\mathcal
L)=(X',\mathcal L')$ and $\L$ is base point free, complete and the
pair $(X,D)$ is as in (i), (ii), (iii) or (iv) of Theorem
\ref{Enr}.\end{Corollary}

\begin{Corollary} Let $X\subset\p^{r}$, $r\geq 3g+5$, $g\geq 2$, be an
irreducible, non--degenerate surface which is not a scroll and
having general hyperplane section $D$ of geometric genus $g$. Then
$r=3g+5$, the surface $X$ is linearly normal, of degree $4g+4$ and
it is one of the following:

\begin{itemize}

\item [(i)] $g=3, r=14$ and $X=V_{2,4}$ is the $4$--Veronese
embedding of $\p^2$ in $\p^{14}$;

\item [(ii)] $X=X_{a,g}$ is a smooth $(a,g)$--Castelnuovo surface,
with $0\leq a\leq g$;\par

\item [(iii)] $X$ has only one singular point and it is the
$2$--Veronese embedding of a cone over a rational normal curve of
degree $a$, $a\geq 3$ and $g=a-1$, i.e. $X=X_{g+1,g}$ is a
$(g+1,g)$-Castelnuovo surface.
\label{Enriquesirreducible}\end{itemize} \end{Corollary}

We finish this section by proving a slight extension of the above
results, which will be essential in our subsequent classification
theorems. Further generalizations, in the spirit of \cite {Ca} or
\cite{Re}, can be obtained, but we will not consider them here,
since we will not use them now. Similarly, we refrain from
formulating the next result in its  maximal generality, i.e. for
big and nef, but not necessarily ample, pairs, since we will not
need such a generality here.

\begin{Theorem} Let $X$ be a smooth, irreducible,
projective surface. Let $D$ be an effective ample divisor on $X$.
Set $d:=D^2$, $g:=p_a(D)$, $r:=\dim(|D|)$. Assume that $g\geq 2$
and that the pair $(X,D)$ is minimal, not a scroll and suppose
that $r=3g+5-s$, with $1\leq s\leq 3$. Then $X$ is rational, $D$
is very ample, and one of the following cases occurs:

\begin{itemize}

\item [(i)] $(X,D)$ is a projection of a 4--Veronese pair from
$i=1,2,3$ points. One has $g=3$, $d=16-s$ and $s=i$;\par

\item [(ii)] $(X,D)$ is a projection of an $(a,g)$--Castelnuovo
pair, with $0\leq a\leq g$, from $i=1,2,3$ points. One has
$d=4g+4-s$, $s=i$;\par

\item [(iii)] $X\simeq \p^1\times \p^1$ and $D$ is of type $(3,3)$
on $X$. One has $g=4$, $d=18$ and $s=2$;\par

\item [(iv)] $(X,D)$ is the tangential projection of a
$5$--Veronese pair from $i=0,1,2$ points. One has $g=6-i$,
$d=25-4i$, $s=3$.\end{itemize}\label{3k+2}
\end{Theorem}

\begin{proof} By the theorem of Riemann--Roch we have $d-g+1\geq r\geq
3g+5-s$, hence $d\geq 4g+4-s$. Moreover, by \eqref{bound}, $d\leq
4g-4+K^2$, so that $K^2\geq 8-s\geq 5$ and $X$ is rational since
$K\cdot D=2g-2-d\leq -2g-1<0$. By \eqref {quadrato}, we have

\begin{equation} \label{degagg} (K+D)^2=K^2-8+s.\end{equation}

Notice that $D^2=d\geq 4g+1\geq 9$ implies, by Reider's theorem
(see \cite{BS}) and the hypotheses  $D$  ample and $(X,D)$ not a
scroll, that $|K+D|$ is base point free. So either $(K+D)^2=0$ and
$|K+D|$ is composite with a base point free pencil $|M|$, or the
general curve $C\in |K+D|$ is smooth and irreducible. Note also
that $\dim(|K+D|)=g-1$. Hence if $g=2$, then $|K+D|$ is a base
point free pencil and therefore $(K+D)^2=0$.\par

Assume that $K^2=9$, i.e. $X\simeq \p^2$. Then (\ref{degagg})
implies that $(K+D)^2=1+s$. So the only possibility is $s=3$ and
$(X,D)$ is a $5$--Veronese pair.\par

From now on we will assume $K^2\leq 8$ and therefore $0\leq
(K+D)^2\leq s\leq 3$ by (\ref{degagg}). We examine separately the
various cases.\par

If $(K+D)^2=0$ and $|K+D|$ is composite with a base point free
pencil $|M|$, the general curve in $|D|$ is hyperelliptic and
therefore $D\cdot M=2$. Since $M\cdot (K+D)=0$, we have $K\cdot
M=-D\cdot M=-2$, and $M^2=0$ yields that the general curve in
$|M|$ is rational. By (\ref {degagg}) we have $K^2=8-s$, so we
have $s$ reducible curves in $|M|$, which are formed by pairs of
$(-1)$--curves meeting transversally at one point and both meeting
$D$ at one point. By contracting $s$ disjoint of these
$(-1)$--curves, we have a morphism $p: X\to \FF_a$, for some
$a\geq 0$. Let $D'=p_*(D)$. Then $p_a(D')=g$ and $D'^2=d+s=4g+4$.
Then, by Theorem \ref{Enr} and Corollary
\ref{Enriquesirreducible}, we conclude we are in case (ii).\par

If $(K+D)^2=1$, then $\phi_{|K+D|}$ is a birational morphism of
$X$ to $\p^2$, hence $X$ is the blow--up of $\p^2$ at $9-K^2=s$
points $x_1,...,x_s$. If $E$ is a $(-1)$--curve contracted by
$|K+D|$, then one has $E\cdot (K+D)=0$, hence $E\cdot D=-E\cdot
K=1$, which means that the image of $|D|$ in $\p^2$ has simple
base points at $x_1,...,x_s$. Furthermore $g-1=\dim(|K+D|)=2$,
hence $g=3$. We are thus in case (i).\par

If $(K+D)^2=2$, then the series cut out by $|K+D|$ on its general
curve $C$ is a complete $g_2^{g-2}$, which implies $g\leq 4$. \par

If $g=4$, then $C$ is rational and $\phi_{|K+D|}$ is a birational
morphism of $X$ to a quadric in $\p^3$. Thus $X$ is the blow--up
of $\FF_a$, $a=0,2$, at $8-K^2=s-2$ points. Note that the
ampleness hypothesis on $D$ rules out the case $a=2$. Then
$s-2\geq 0$, namely $2\leq s\leq 3$. If $s=2$, then we clearly are
in case (iii), whereas, if $s=3$, we are in case (iv), $i=2$. \par

Suppose $g=3$. Let $C$ be the general curve in $|K+D|$. One
computes $(K+C)\cdot C=0$ and $(K+C)^2=(2K+D)^2=8-s>0$. This
contradicts the Hodge index theorem.\par

If $(K+D)^2=3$, then the series cut out by $|K+D|$ on its general
curve $C$ is a complete $g_3^{g-2}$, which implies $g\leq 5$. On
the other hand (\ref{degagg}) implies that $s=3$, $K^2=8$, i.e.
$X$ is a surface $\FF_a$, for some $a\geq 0$.\par

If $g=5$, then $C$ is rational and $\phi_{|K+D|}$ is then an
isomorphism of $X$ to $\FF_1$ embedded in $\p^4$ as a rational
normal cubic scroll. It is then clear that we are in case (iv),
$i=1$. \par

If $g\leq 4$, one computes $(K+C)\cdot C=8-2g$ and
$(K+C)^2=21-4g$, which contradicts the Hodge index theorem.\par

The proof is thus completed.\end{proof}

The pairs listed in (i)--(iv) of Theorem \ref {3k+2} above will be
called {\it almost extremal Castelnuovo pairs}. The corresponding
surfaces $\phi_{|D|}(X)$ will be called {\it almost extremal
Castelnuovo  surfaces}.

\section{The classification of $\mathcal{OA}^{k+1}_{k-1}$--surfaces}
\label{oasurf}

In this section we give the classification of surfaces
$X\subset\p^{3k+2}$, $k\geq 2$ with $\nu_k(X)=1$. Recall that the
case $k=1$ was classically considered by Severi \cite {Se} and
proved by the second author \cite {Ru} (see also \cite {CMR}). We
notice that this classification was in part divined by Bronowski
in \cite{Br}, where however the argument he gives relies on the
unproved conjecture stated in Remark \ref {Bronowski}.

\begin{Theorem}  Let $X\subset\p^{3k+2}$,
$k\geq 2$, be a smooth, projective, surface which is linearly
normal, and such that $\nu_k(X)=1$. We let $d$ be the degree and
$g$ be the sectional genus of $X$. Then $X$ is one of the
following:

\begin{itemize}

\item [(i)] a rational normal scroll $S(a_1,a_2)$  with $k\leq
a_1\leq a_2$, $d=a_1+a_2=3k+1$ and sectional genus $g=0$ (see
Example \ref {rationalscrolls});\par

\item [(ii)] an odd Castelnuovo surface $X_{a,k-1}$, with $0\leq
a\leq k-1$ and $a+k\equiv 1$ (mod. 2) (see Example
\ref{quadricfibrationsextended} and Remark \ref {castsurf}). In
this case $d=4k$, $g=k-1$ and the hyperplane sections of $X$ are
hyperelliptic curves;\par

\item [(iii)] the internal projection from three distinct points
of a Castelnuovo surface $X_{a,k} \subset\p^{3k+5}$ with $0\leq
a\leq k$. In this case $d=4k+1$ and $g=k$ and the hyperplane
sections are hyperelliptic curves (see Example
\ref{quadricfibrationsextended});\par

\item [(iv)] the tangential projection of a $5$--Veronese surface
$V_{2,5}$ from $i=0,1,2, 3$ points (see Example \ref {quintics}).
Here $d=25-4i$, $g=k=6-i$.\end{itemize}\label{nu=1} \end{Theorem}

\begin{proof} From the classification of weakly defective surfaces (see
\cite {CC}, Theorem 1.3 and Example \ref {wdefectivesurf} above),
we see that $X$, being not $k$--defective and spanning a
$\p^{3k+2}$, is also not $k$--weakly defective. We can, and will,
therefore apply Proposition \ref {tangent}. Let $p_1,...,p_k\in X$
be general points and let $\mathcal L$ be the linear system of
hyperplane sections of $X$ tangent at $p_1,...,p_k$. Since $X$ is
not $(k-1)$--defective, we have $\dim({\mathcal L})=2$. Moreover
${\mathcal L}=F+{\mathcal M}$, where $F$ is the fixed part and
$\mathcal M$ the movable part, as described in Proposition \ref
{tangent}. The relevant information is that, by Theorem \ref{boh},
$\tau_{X,k}: X\map \p^2$ is birational, hence $X$ is rational and
the general curve $M\in {\mathcal M}$ is rational and $\mathcal M$
determines a birational map of $X$ to $\p^2$. In particular,
${\mathcal M}$ is base point free off $p_1,...,p_k$ (see
\cite{CMR}, Proposition 6.3). \par

We will separately discuss the various cases according to
Proposition \ref {tangent}:

\begin{enumerate}

\item $F$ is empty;

\item $F$ is not empty and irreducible;

\item $F$ consists of $k$ irreducible curves $\Gamma_i$ with
$p_i\in \Gamma_i$.

\end{enumerate}

In case (1) the curve $M$ is rational with $k$ nodes at
$p_1,...,p_k$ and no other singularity. Then $g=k$ and $d=4k+1$
and therefore $X$ is an almost extremal Castelnuovo surface with
$\epsilon=3$. By Theorem \ref {3k+2}, we are either in case (iii)
or in case (iv).\par

In case (2), the curve $F$ is smooth and rational. Look at the
linear system $|F|$ on $X$. Since $X$ is linearly normal and there
is a unique curve $F$ containing the general points $p_1,...,p_k$,
then we have $\dim(|F|)=k$, hence $F^2=k-1$. Moreover $M$ is also
rational and smooth. Look at the system $|M|$. Since there is a
$2$--dimensional linear system of curves in $|M|$ containing
$p_1,...,p_k$, we have $\dim(|M|)=k+2$, thus $M^2=k+1$. Moreover
$M\cdot F=k$ by  Proposition \ref {tangent}. This implies that:

$$d=M^2+2M\cdot F+F^2=4k,\quad g=p_a(M)+p_a(F)+M\cdot F-1=k-1$$

\noindent hence $X$ is an extremal Castelnuovo surface. By
Corollary \ref {Enriquesirreducible}, we are in case (ii), because
the Veronese surface $V_{2,4}$ is $4$--defective (see Remark
\ref{castsurf}). \par

In case (3), the curves $\Gamma_i$ are rational and linearly
equivalent, and $\Gamma_i^2=0$, for $i=1,...,k$. This
 implies that
we are in case (i).\end{proof}

\begin{Remark} {\rm The assumption that $X$ be linearly
normal is essential to have a {\it finite} classification in
Theorem \ref {nu=1} above, as shown in Example \ref {wdd}. We do
not know whether there are more examples of non--linearly normal
$\mathcal{OA}^{k+1}_{k-1}$--surfaces other than the ones exhibited
in Example \ref {wdd}.

According to Proposition \ref{linearnorm}, $k$--regularity implies
linear normality. So one could be tempted to replace the linear
normality hypothesis in Theorem \ref {nu=1} by the $k$--regularity
assumption, which seems to be, in this context, a right
generalization of the concept of smoothness. However the
$k$--regularity hypothesis is almost never verified by the
surfaces in the list (i)--(iv) of Theorem \ref {nu=1}. This
suggests that $k$--regularity is too rigid. It would be
interesting to find a weaker concept which, in this context, could
play the right role.}\label{kreg} \end{Remark}

\section{The classification of ${\mathcal M}^k$--surfaces}\label{msurf}

In this section we consider the classification of ${\mathcal
M}^k$--surfaces (see also \cite {Br}). The case of $k$--defective
and $k$--weakly defective surfaces has been already considered in
Examples \ref {defectivesurf}, \ref {wdefectivesurf} and
\ref{wdd}. We summarize the result in the following:

\begin{Theorem} \label{kdefect} Let $X\subset \p^r$ be
an irreducible, non--degenerate, surface. If $X$ is
$k$--defective, then it is an ${\mathcal M}^k$--surface if and
only if one if the following happens:\par

\begin{itemize}

\item [(i)]  $r=3k+2$ and $X$ is the $2$--Veronese embedding of a
surface of degree $k$ in $\p^{k+1}$;\par

\item [(ii)] $X$ sits in a $(k+1)$--dimensional cone, with a
vertex of dimension $k-1$, over a rational normal curve $C$ of
degree $d\geq 2k+3$. \par

\end{itemize}

If $X$ is $k$--weakly defective, but not $k$--defective, then it
is an ${\mathcal M}^k$--surface if and only if one if the
following happens:\par

\begin{itemize}

\item [(iii)] $r=9$, $k=2$ and $X$ is the $2$--Veronese embedding
of a surface of degree $d\geq 3$ in $\p^{3}$;\par

\item [(iv)]  $r=3k+3$ and $X$ is the cone over a $k$--defective
surface of type (i);\par

\item [(v)] $r=3k+3$ and $X$ is the $2$--Veronese embedding of a
surface of degree $k+1$ in $\p^{k+1}$;\par

\item [(iv)]  $X$ sits in a $(k+2)$--dimensional cone, with a
vertex of dimension $k$, over a rational normal curve $C$ of
degree $d\geq 2k+2$.
\end{itemize}\end {Theorem}

The main result of this section is the classification theorem for
$\mathcal{MA}^{k+1}_{k-1}$--surfaces, which concludes the
classification of ${\mathcal M}^k$--surfaces:

\begin{Theorem} \label{kdefects} Let $X\subset\p^{3k+2+h}$, with $k,h\geq
1$, be a smooth, irreducible, non--degenerate,
$\mathcal{M}^{k}$--surface which is linearly normal and not
$k$--weakly defective. Let $d$ be the degree and $g$ the sectional
genus of $X$. Then  $X$ is  one of the following:\par

\begin{itemize}

\item [(i)] a rational normal scroll $S(a_1,a_2)$ of degree
$d=3k+1+h$ and type $(a_1,a_2)$ with $k\leq a_1\leq a_2$ (see
Example \ref {rationalscrolls});  \par

\item [(ii)] a del Pezzo surface of degree $d=5+h$  and $g=1$,
with $1\leq h\leq 4$ and $k=1$ (see Example \ref {delpezzo});\par

\item [(iii)] the internal projection from $3-h$, with $1\leq
h\leq 3$, distinct points  of an odd Castelnuovo surface $X_{a,k}
\subset\p^{3k+5}$ with $0\leq a\leq k$ and $a+k\equiv 0$ (mod. 2).
In this case $d=4k+1+h$, $g=k$ and the hyperplane sections are
hyperelliptic curves (see Example
\ref{quadricfibrationsextended});\par

\item [(iv)] the internal projection from $3-h$ points, with
$1\leq h\leq 2$, of the Veronese surface $V_{2,4}$. In this case
$d=13+h, g=3, k=3$ (see Example \ref {quartics});\par

\item [(v)] the $3$-Veronese embedding in $\p^{15}$ of a smooth
quadric in $\p^3$. Here $d=18, g=4, k=4, h=1$ (see Example
\ref{treverquadrics});\par

\item [(vi)] the 3-Veronese embedding $V_{2,3}$ of $\p^2$. In this
case $d=9$, $g=1$, $k=2, h=1$ (see Example \ref {delpezzo}).
\end{itemize}\label{mak}
\end{Theorem}

\begin{proof}  Since $X$ is not $k$--weakly defective, we can apply again
Proposition \ref {tangent}. Let $p_1,...,p_k\in X$ be general
points and, as in the proof of Theorem \ref {nu=1}, we let
$\mathcal L$ be the linear system of hyperplane sections of $X$
tangent at $p_1,...,p_k$. Since $X$ is not $(k-1)$--defective, we
have $\dim({\mathcal L})=2+h$. Moreover ${\mathcal L}=F+{\mathcal
M}$, where $F$ is the fixed part and $\mathcal M$ the movable
part, as described in Proposition \ref {tangent}. By  Corollary
\ref{minimal}, $\tau_{X,k}: X\map X_k\subset \p^{h+2}$ is
birational and $X_k$ is a surface of minimal degree $h+1$, hence
$X$ is rational and the general curve $M\in {\mathcal M}$ is also
rational. \par

Again, as in the proof of Theorem \ref {nu=1}, one has to
separately discuss the various cases according to Proposition \ref
{tangent}.

If $F$ is empty, then $g=k$ and $d=4k+h+1$. If $k=1$ we are in
case (ii). If $k>1$, by applying  Corollary
\ref{Enriquesirreducible} and Theorem \ref {3k+2}, we see that we
have cases (iii), (iv) and (v).\par

If $F$ is not empty and irreducible, then $g=k-1$ and $d=4k+h$. By
Theorem \ref {Enr}, the only possible case is $h=1$, $g=1$, which
implies $k=2$ and we are in case (vi).

If $F$ consists of $k$ irreducible curves we are in case
(i).\end{proof}

We con now state our result concerning the generalized Bronowski's
conjecture for surfaces (see Remark \ref {Bronowski}):

\begin{Corollary} \label{Bron} The generalized Bronowsi's conjecture
holds for smooth surfaces.
\end{Corollary}

\begin{proof} Let $X\subset \p^{3k+2+h}$, $h:=h^{(k)}(X)$, be a smooth,
irreducible, projective, not $k$--defective surface and assume
that the general $k$--tangential projection $\tau_{X,k}: X\map
X_k\subset \p^{h+2}$ birationally maps $X$ to a surface of minimal
degree $h+1$ in $\p^{h+2}$. The same argument we made in the
proofs of Theorems \ref {nu=1} and \ref {kdefect} proves that $X$
is either or minimal degree or Castelnuovo extremal or Castelnuovo
almost extremal. As we saw in \S \ref {ex}, these are
$\mathcal{MA}^{k+1}_{k-1}$ or
$\mathcal{OA}^{k+1}_{k-1}$--surfaces, according to whether $h>0$
or $h=0$.
\end{proof}

\section{A generalization of a theorem of Severi}\label{gensev}

Terracini's Lemma \ref {terracini} implies that a defective
variety is swept out by very degenerate subvarieties. As a
consequence, one has a famous theorem of Severi \cite{Se} (see
also \cite{Ru2}), which says that the Veronese surface $V_{2,2}$
in $\p^5$ is the only irreducible non--degenerate, projective
surface in $\p^r$, $r\geq 5$, not a cone, such that
$\dim(S(X))=4$. This result can be restated  as follows: {\it the
Veronese surface in $\p^5$ is the only 1--defective, not 0--weakly
defective, irreducible non--degenerate, projective surface in
$\p^r$, $r\geq 5$} (cfr. Remark \ref{0wd}).\par

This section is devoted to point out an extension of Severi's
theorem, namely Theorem \ref {Bronowskidifettivo} below. This
result yields a projective characterization of extremal
Castelnuovo surfaces, in particular it stresses a distinction
between odd and even $(a,k)$-Castelnuovo surfaces, as suggested by
Bronowski in \cite{Br}.\par

Theorem \ref{Bronowskidifettivo} could also be deduced by the
classification of weakly defective surfaces (see \cite {CC} and
Examples \ref {defectivesurf} and \ref {wdefectivesurf}). However
the proof in \cite {CC} requires a subtle analysis involving
involutions on irreducible varieties and a generalization of the
Castelnuovo--Humbert theorem to higher dimensional varieties. It
seems interesting to us to present here an easy argument based on
the ideas developed in this paper.

\begin{Theorem}\label{Bronowskidifettivo}
Let $X\subset\p^r$, $r\geq 3k+2$ and $k\geq 1$, be a smooth,
irreducible, non--degenerate surface. Suppose that
 $X$ is $k$--defective but not $(k-1)$--weakly defective.
 Then $r=3k+2$ and
 $X$ is the $2$--Veronese embedding of a smooth surface of degree $k$ in
 $\p^{k+1}$, i.e.  it is one of the following:

\begin{itemize}
\item[(i)] $X=V_{2,2}$ is the Veronese surface in $\p^5$, then
$k=1$ and $\deg(S(X))=3$;

\item [(ii)]  $X=V_{2,4}$ is the $4$--Veronese embedding of $\p^2$
in $\p^{14}$, then $k=4$ and $\deg(S^4(X))=6$;

\item [(iii)] $X$ is a smooth even Castelnuovo surface
$X_{a,k-1}$, with $0\leq a\leq k-1$, which is the $2$--Veronese
embedding of a smooth rational normal scroll of degree $k$ in
$\p^k$.\end{itemize}

In particular a $k$-defective, not $(k-1)$--weakly defective,
surface in $\p^r$, $r\geq 3k+2$, is an $\mathcal M^k$--surface in
$\p^{3k+2}$.

\end{Theorem}

\begin{proof} Let $p_0,...,p_k\in X$ be general points. Since $X$ is not
$(k-1)$--defective, one has $\dim(T_{X,p_1,...,p_k})=3k-1$. Since
$X$ is not degenerate in $\p^r$, $r\geq 3k+2$, the projection of
$X$ from $T_{X,p_1,...,p_k}$ cannot be a point. Hence
$s^{(k)}(X)=\dim(T_{X,p_0,...,p_k})=3k+1$.\par

We can suppose $k\geq 2$ by Severi' theorem  \cite{Se}. Also we
may assume that $X\subset\p^r$ is linearly normal. Since $X$ is
not $(k-1)$--weakly defective we may apply Lemma \ref{birational}
to deduce that $\tau_{X,k-1}: X\map X_{k-1}\subset \p^{r-3k+3}$ is
birational to its image. Then $r-3k+3\geq 5$ and
$X_{k-1}\subset\p^{r-3k+3}$ is an irreducible non--degenerate
surface. By Terracini's Lemma $\dim(S(X_{k-1}))=4$ and moreover
$X_{k-1}$ is not $0$--weakly defective because $X\subset\p^r$ is
not $(k-1)$--weakly defective. Thus Severi's theorem applies and
yields that $X_{k-1}$ is the Veronese surface in $\p^5$ and that
$r=3k+2$. Note that $X$ cannot be a scroll, since
$X_{k-1}=V_{2,2}$ does not contain lines.\par

The rest of the proof is analogous to the one in Theorem
\ref{nu=1}. Since $X$ is not $(k-2)$--weakly defective we can
apply Proposition \ref{tangent}. Let $p_1,...,p_{k-1}\in X$ be
general points and, as in the proof of Theorem \ref{nu=1}, we let
$\mathcal L$ be the linear system of hyperplane sections of $X$
tangent at $p_1,...,p_{k-1}$. The general curve $M\in {\mathcal
M}$ is rational being birational to a hyperplane section of the
Veronese surface $X_{k-1}\subset\p^5$ and we have $\dim({\mathcal
L})=5$. Moreover ${\mathcal L}=F+{\mathcal M}$, where $F$ is the
fixed part and $\mathcal M$ the movable part, as described in
Proposition \ref {tangent}.
\par

Again, one has to separately discuss the various cases according
to Proposition \ref {tangent}.

If $F$ is empty, then $g=k-1$ and $d=4k$. In the case $k=2$, then
$X$ is a del Pezzo surface of degree $8$ and we are in case (iii)
(see Example \ref {delpezzo}). If $k\geq 3$, by applying Corollary
\ref {Enriquesirreducible}, we have cases (ii) and (iii).\par

If $F$ is not empty and irreducible, then $g=k-2$ and $d=4k-1$. We
can suppose that $k\geq 3$ since $X$ is not a scroll. Note also
that $3(k-2)+5=3k-1$. Since $X$ is not a scroll, then Corollary
\ref{Enriquesirreducible} implies that this case does not exist.

If $F$ consists of $k-1$ irreducible curves, then they belong to a
 pencil of lines,  a contradiction, since $X$ is not a scroll.\end{proof}

\end{document}